\documentclass[11pt,twoside]{article}

\AtBeginDocument{%
  \fontsize{11pt}{13pt}\selectfont
}

\makeatother

\renewenvironment{abstract}
  {\small
   \begin{quote}
   \noindent\textbf{Abstract.}\ }
  {\end{quote}}
\newcommand{\keywords}[1]{%
  \par\medskip
  \noindent\textit{Keywords:} #1\par
}

\newcommand{\subjclass}[1]{%
  \smallskip
  \noindent\textit{2010 Mathematics Subject Classification:} #1\par
}
\usepackage{amsthm}

\theoremstyle{plain}
\newtheorem{lemma}{Lemma}[section]
\newtheorem{theorem}[lemma]{Theorem}

\theoremstyle{definition}

\theoremstyle{remark}
\newtheorem{remark}[lemma]{Remark}

\usepackage{graphicx}
\usepackage{titlesec}
\usepackage{newtxtext,newtxmath}

\titleformat{\section}[block]
  {\centering\large\bfseries\color{black}}
  {\thesection.}{0.6em}{}

\titleformat{\subsection}[hang]
  {\normalfont\large\bfseries\color{black}}
  {\thesubsection}{1em}{}

\titleformat{\subsubsection}[hang]
  {\normalfont\normalsize\bfseries\color{black}}
  {\thesubsubsection}{1em}{}

\usepackage[T1]{fontenc}
\usepackage{amsmath}
\usepackage[margin=0.87in]{geometry}
\usepackage{bm}
\usepackage{xcolor}

\definecolor{A}{RGB}{138,43,226}
\definecolor{C}{RGB}{0,0,254}
\definecolor{B}{RGB}{25,25,112}
\definecolor{S}{RGB}{20,0,240}
\definecolor{G}{RGB}{70,130,180}
\definecolor{R1}{RGB}{139,58,58}
\definecolor{R2}{RGB}{208,32,144}
\definecolor{R3}{RGB}{0,139,139}
\definecolor{br}{RGB}{135,38,87}

\usepackage[
  breaklinks=true,
  colorlinks=true,
  linkcolor=C,
  urlcolor=blue,
  citecolor=C,
  bookmarks=true,
  bookmarksopen=true,
  bookmarksdepth=3
]{hyperref}
\usepackage{orcidlink}
\usepackage[
  backend=biber,
  style=numeric-comp,
  sorting=nyt,
  sortcites=true
]{biblatex}

\usepackage{fancyhdr}

\fancypagestyle{plain}{%
  \fancyhf{}

}

\usepackage{xcolor}

\usepackage[backend=biber]{biblatex}
\usepackage{float}
\usepackage{tocloft}

\title{\textbf{PV-Wave Decomposition and Combined Wave Effects for Stratified Couette Flows}}
\author{Yiting Yao\,\orcidlink{0009-0000-5232-7230}}
\date{}

\makeatletter
\renewcommand{\@maketitle}{%
  \newpage
  \vspace*{3pt}
  \noindent{\LARGE\@title\par}
  \vskip 1em
  \noindent{\large\bfseries\@author\par}
  \vskip 0.4em
  \noindent{\small
    Department of Mathematics, University of Bath, Bath, UK\\
    E-mail: yy2644@bath.ac.uk\par
  }
  \vspace{-1pt}
}
\makeatother

\begin{document}

\maketitle
\vspace{1pt}

\begin{abstract}
Motivated by a class of three-component Fourier ODE systems consisting of a conserved mode coupled to a pair of oscillatory modes, we develop a "potential-vorticity-wave" (PV-Wave) decomposition strategy for constructing adapted energy variables in the stability analysis of Couette flows. The strategy separates the conserved potential-vorticity component, symmetrizes the remaining wave subsystem, and introduces an additional correction to control the couplings generated by time-dependent Fourier coefficients. We apply this framework to study the linearized dynamics of the 3D Boussinesq MHD system and rotating Boussinesq systems for Couette flows and quantify the inviscid damping and parameter-dependent amplification arising from the interaction of Alfvén-gravity waves and inertial-gravity waves.

\end{abstract}

\section{Introduction}
The interaction of magnetic fields, rotation, and buoyancy is ubiquitous in astrophysical and geophysical fluid dynamics, for instance in the evolution of Earth's core and the propagation of internal-gravity waves in the atmosphere. A standard model for these phenomena is the rotating Boussinesq MHD system, given by\begin{flalign}
\left\{
\begin{aligned}
& \partial_t v+(v\cdot\bm{\nabla})v+\widetilde{\gamma}\vec{\bm{e}}_3\times v=-\bm{\nabla}p+(B\cdot\bm{\nabla})B+\varrho\vec{\bm{e}}_3,\\
& \partial_t\varrho+v\cdot\bm{\nabla}\varrho=0,\\
& \partial_t B+(v\cdot\bm{\nabla})B-(B\cdot\bm{\nabla})v=0,\\
& \bm{\nabla}\cdot v=0,\qquad \bm{\nabla}\cdot B=0.
\end{aligned}
\right.
\label{RBM}
\end{flalign}Here $v,B\in\mathbb{R}^3$ denote the fluid velocity and magnetic field, respectively. Both vector fields are divergence free, with $\bm{\nabla}\cdot v=0$ expressing the incompressibility of the fluid. The scalar fields $\varrho$ and $p$ denote the buoyancy variable and pressure. $\widetilde{\gamma}\in\mathbb{R}$ is the Coriolis parameter. In the absence of a magnetic field, that is, when $B\equiv0$, system \eqref{RBM} admits a materially conserved quantity\footnote{A quantity $s$ (scalar or vector) is materially conserved if $\partial_t s+v\cdot\bm{\nabla}s=0.$}, known as (Ertel's) \textit{potential vorticity} (PV):\begin{equation}
q_1:=(\bm{\omega}+\widetilde{\gamma}\vec{\bm{e}}_3)\cdot\bm{\nabla}\varrho,\label{potentialvorticity}
\end{equation}where $\bm{\omega}:=\bm{\nabla}\times v\in\mathbb{R}^3$ is the vorticity of the fluid. PV is a key quantity in rotating stratified fluids arising in the atmosphere and ocean. It has important applications in geophysical fluid dynamics, such as weather prediction, and can be used to describe the formation and evolution of Rossby waves, jet streams, cyclones, and blocking events \cite{Vallis_2017}. In the presence of magnetic field/Lorentz force, Ertel's PV is no longer materially conserved, and there is a magnetic analogue of potential vorticity, which we call the \textit{magnetic potential vorticity}  \cite{10.1111/j.1365-246X.1996.tb06529.x,PhysRevE.85.026301}
\begin{equation}
q_2:=B\cdot\bm{\nabla}\varrho.\label{magneticpv}
\end{equation}
The subject of this work is to investigate how these PV invariants can be used in the stability analysis of the steady states. The stability theory of steady states such as shear flow or vortices of fluids has been studied extensively, with particular attention to the dynamics near Couette flow. In their pioneering work, Bedrossian and Masmoudi \cite{bedrossian2015inviscid} proved the nonlinear inviscid damping of the Couette flow in 2D Euler. Similar developments had then been made on the stability of Couette flow for 3D Navier Stokes by Bedrossian, Masmoudi and Germain \cite{bedrossian2015dynamicsnearsubcriticaltransition,bedrossian2015dynamicsnearsubcriticaltransition1,9e0bb666-1e4e-3aa2-87b0-f916c2336491}. There have been several works later on, establishing the stability of the Couette flow in different models such as Boussinesq, (in)compressible Navier Stokes, MHD, see, for instance \cite{li2025transitionthresholdnavierstokescoriolishigh,bedrossian2018sobolev,nonlinearinvisciddamping,coti2024stability,li2026nonlinearstabilitythreshold3d,KNOBEL2026113937,liss2020sobolev,antonelli2021linear,yao2026stabilitystratifiedcouette}. In the companion work \cite{yao2026stabilitystratifiedcouette},
we study the non-parallel
configuration consisting of the plane Couette flow
$v^*=(y,0,0)$ and vertical stratification. There, the linearized PV
is used to reformulate the non-zero-mode system, reducing the
original three-component system to a symmetrizable two-component
system, with the PV acting as an explicitly controlled forcing. This PV-based reformulation is also a key ingredient
in the proof of the nonlinear stability of $v^*$.\vspace{9pt}

In this article, we are concerned with the following two problems. \begin{enumerate}
\renewcommand{\labelenumi}{(\arabic{enumi})}
    \item (\textit{Boussinesq MHD}) In the presence of magnetic field and absence of rotation (i.e. $\widetilde{\gamma}=0$ in \eqref{RBM}), we consider the following equilibrium\begin{equation}
    v^*=(z,0,0),\qquad B^*=(0,\beta,0),\qquad \varrho^*=\alpha z.\label{state1}
    \end{equation}The velocity field is the vertical Couette flow, and the
    background magnetic field points in the $y$-direction. Stable stratification is characterized by $\alpha>0$, with \textit{Brunt-V\"ais\"al\"a} frequency $\sqrt{\alpha}$. We will consider perturbations on $\mathbb{T}^2\times\mathbb{R}$. 
    \item (\textit{Rotating Boussinesq}) In the absence of magnetic field ($B\equiv 0$ in \eqref{RBM}), we study the following equilibrium\begin{equation}
    v^*=(y,0,0),\qquad \varrho^*=\alpha z,\qquad \alpha>0\label{state2}
    \end{equation}Here $v^*$ is known as the plane Couette flow, and we consider perturbations on $\mathbb{T}\times\mathbb{R}\times\mathbb{T}$.
\end{enumerate}\vspace{5pt}

\noindent In both cases, we write the perturbations for $v$ and $\varrho$ as\[
v=v^*+u,\qquad \varrho=\varrho^*+\sqrt{\alpha}\theta. 
\]In the Boussinesq MHD case, we also write $B=B^*+b$. The linearized equations around \eqref{state1} are\begin{flalign}
\left\{
\begin{aligned}
& \partial_t u+z\partial_x u+u_3\vec{\bm{e}}_1+\bm{\nabla}p=\beta\partial_yb+\sqrt{\alpha}\theta\vec{\bm{e}}_3,\\
& \partial_t\theta+z\partial_x\theta+\sqrt{\alpha}u_3=0,\\
& \partial_t b+z\partial_x b=\beta\partial_yu+b_3\vec{\bm{e}}_1,\\
& \bm{\nabla}\cdot u=0,\qquad \bm{\nabla}\cdot b=0.
\end{aligned}
\right.
\label{LBM}
\tag{LBM}
\end{flalign}Similarly, the linearized system around state \eqref{state2} is given by\begin{flalign}
\left\{
\begin{aligned}
&\partial_t u+y\partial_x u+u_2\vec{\bm{e}}_1+\widetilde{\gamma}\vec{\bm{e}}_3\times u+\bm{\nabla}p=\sqrt{\alpha}\theta\vec{\bm{e}}_3,\\
    &\partial_t\theta+y\partial_x\theta+\sqrt{\alpha} u_3=0,\\
     &\bm{\nabla}\cdot u=0.
\end{aligned}
\right.
\label{LRB}
\tag{LRB}
\end{flalign}Our aim is to establish quantitative bounds for solutions to \eqref{LBM} and \eqref{LRB}, thereby characterizing inviscid damping, dispersive decay, and instability in the relevant parameter regimes.

Throughout linear stability analysis, the main difficulty arises in the analysis of the coupled system $(u_3,\theta,b_3)$ in the first problem and $(u_3,\theta,\omega_3)$ in the second. To remove the background advection, we introduce the following moving coordinates: for the Boussinesq MHD problem, we define\[
X=x-tz,\qquad Y=y, \qquad Z=z.
\]For the rotating Boussinesq system, we set\[
X=x-ty,\qquad Y=y,\qquad Z=z.
\]Taking the Fourier transform in $(X,Y,Z)$, with corresponding frequency variables $(k,\eta,l)$, we derive the $(\widehat{U}_3,\widehat{\Theta},\widehat{B}_3)$ system\footnote{We use capitalized letters to denote the function in the moving frame, e.g. $M(t,X,Y,Z)=m(t,x,y,z)$.}\begin{flalign}
\left\{
\begin{aligned}
&\partial_t \widehat{U}_3=-\frac{p'}{p}\widehat{U}_3+i\beta \eta \widehat{B}_3+\sqrt{\alpha}\frac{f}{p}\widehat{\Theta},\\
&\partial_t\widehat{\Theta}=-\sqrt{\alpha}\widehat{U}_3,\\
&\partial_t \widehat{B}_3=i\beta \eta\widehat{U}_3,
\end{aligned}
\right.\label{system1}
\end{flalign}and the $(\widehat{U}_3,\widehat{\Theta},\widehat{W}_3)$ system\begin{flalign}
    \left\{
\begin{aligned}
&\partial_t  \widehat{U}_3=\frac{p'l^2}{pf}\widehat{U}_3+\frac{i\widetilde{\gamma}l}{p}\widehat{W}_3+\sqrt{\alpha}\frac{f}{p}\widehat{\Theta}-\frac{2ik^2l}{pf}\widehat{W}_3,\\
&\partial_t\widehat{W}_3=-il(1-\widetilde{\gamma})\widehat{U}_3,\\
    &\partial_t\widehat{\Theta}=-\sqrt{\alpha}\widehat{U}_3.
\end{aligned}
\right.\label{system2}
\end{flalign}

Here $f$ and $p$ are the Fourier symbols for the negative horizontal and usual Laplacian in the moving frame, i.e. $-\Delta_H$ and $-\Delta_L$ (see Section  \ref{notation} for the definitions). $p'=\partial_t p$ is the time derivative of $p$. The unknown $W_3$ is understood as the perturbed vertical vorticity in the moving frame. A standard approach to study the decay/boundedness properties is to use the symmetrization method (see, e.g. \cite{bianchini2022linear,antonelli2021linear}). In the present problems, direct symmetrization is ineffective because of the non-integrable-in-time coefficient $\frac{p'}{p}$; estimating this term directly may produce spurious growth in time or an artificial loss of regularity. To resolve this issue, we need to introduce suitable reformulation. We begin with the following observation: consider the Boussinesq MHD problem, the background magnetic field $B^*$ can produce Alfv\'en waves, while stable stratification generates internal-gravity waves. On the other hand, the linearization of the magnetic PV around steady state \eqref{state1} is given by\[
q_{2,l}:=B^*\cdot\bm{\nabla}(\sqrt{\alpha}\theta)+b\cdot\bm{\nabla}\varrho^*.
\]This quantity is conserved along the Couette flow, and it is "non-oscillatory". A natural idea is then to split the system into oscillatory (wave) part and non-oscillatory (PV) part and see how it helps with our analysis.

\subsection*{Motivation for the PV-Wave decomposition.}
To illustrate our idea, we fix $\phi\in(0,\pi/2)$ and consider a simpler toy model:\begin{flalign}
\left\{
\begin{aligned}
 \partial_t U_3&=-\cos\phi\Theta+\sin\phi N,\\
    \partial_t \Theta&=\cos\phi U_3,\\
    \partial_t N&=-\sin\phi U_3.
\end{aligned}
\right.\label{samplesystem}
\end{flalign} Clearly, systems in \eqref{system1} and \eqref{system2} have similar structure to this toy model. Solutions to system \eqref{samplesystem} satisfy conserves the energy\begin{equation}
|U_3|^2+|\Theta|^2+|N|^2=|U_3(0)|^2+|\Theta(0)|^2+|N(0)|^2.\label{sample}
\end{equation}Moreover, $U_3$ solves the following equation\begin{equation}
\partial^2_t U_3=-U_3\qquad \Rightarrow\qquad U_3=a\cos t+b\sin t.\label{secondorder}
\end{equation}Although this system can be solved explicitly, suppose, for the purpose of motivating the decomposition, that neither the explicit solutions nor the energy estimate such as \eqref{sample} can be easily derived, and we only know that the general solution consists of oscillatory (corresponds to the linear combination of $\cos t$ and $\sin t$ here) and non-oscillatory parts (corresponds to the constant part in the solution to \eqref{samplesystem}). The aim is to derive an alternative form of energy estimate. We begin by noticing the existence of a "constant" part:\[
\boxed{\text{PV}:=\sin\phi\Theta+\cos\phi N},\qquad  \Rightarrow\qquad \partial_t \text{PV}=0.
\]We refer to this quantity as the "PV" component of the system. Besides, we can also construct a new unknown "Wave", defined as\[
\boxed{\text{Wave}=A \Theta+B N},
\]where $A$ and $B$ are some constants to be determined. The motivation for defining this unknown is to capture the oscillatory behavior of $\Theta$ and $N$. To make sense of the unknown "Wave", $A$ and $B$ should be chosen such that PV and Wave are linearly independent, i.e.\[
\det\begin{pmatrix}
    \sin\phi&\cos\phi\\
    A& B
\end{pmatrix}=B\sin\phi-A\cos\phi\neq0.
\]A simple choice is then $B=\sin\phi$ and $A=-\cos\phi$ so that we can express $\Theta$ and $N$ in terms of PV and Wave: \begin{equation}
\Theta=\sin\phi\cdot\text{PV}-\cos\phi\cdot\text{Wave},\qquad N=\cos\phi\cdot\text{PV}+\sin\phi\cdot\text{Wave}.\label{pwd2}
\end{equation}We can also derive a closed system in terms of PV, Wave and $U_3$:\begin{flalign*}
\left\{
\begin{aligned}
\partial_t \text{Wave}&=-\omega U_3,\\
    \partial_t \text{PV}&=0,\\
    \partial_t U_3&=\omega  \text{Wave}.
\end{aligned}
\right.
\end{flalign*}Here $\omega:=\sqrt{\cos^2\phi+\sin^2\phi}=1$ can be interpreted as the "wave frequency", it can be obtained from equation \eqref{secondorder}. Moreover, we can derive an alternative form of energy estimate\[
\partial_t\left(|\text{Wave}|^2+|\text{PV}|^2+|U_3|^2\right)=0.
\]From the constructions, we can easily obtain the "norm equivalence"\[
|\text{PV}|+|\text{Wave}|\lesssim|\Theta|+|N|\lesssim |\text{PV}|+|\text{Wave}|,
\]the estimate \eqref{sample} can then be recovered with "$=$" being replaced by "$\lesssim$".\vspace{5pt}

We expect the underlying decomposition principle to extend to more complicated systems by identifying components with qualitatively distinct dynamics. To illustrate this point, consider
\begin{equation*}
\left\{
\begin{aligned}
\partial_t U_3
    &=\mathcal{L}_1U_3-\cos\phi\,\Theta+\sin\phi\,N,\\
\partial_t\Theta
    &=\cos\phi\,U_3+\mathcal{L}_2\Theta,\\
\partial_tN
    &=-\sin\phi\,U_3+\mathcal{L}_3N,
\end{aligned}
\right.
\end{equation*}
where $\mathcal{L}_i$ are time-independent linear operators acting on
the spatial variables. Suppose that
$\mathcal{L}_2=\mathcal{L}_3=\mathcal{L}$, and introduce
\[
P:=\sin\phi\,\Theta+\cos\phi\,N,
\qquad
W:=-\cos\phi\,\Theta+\sin\phi\,N.
\]
A direct computation gives
\[
\partial_tP=\mathcal{L}P,\qquad
\partial_tU_3=\mathcal{L}_1U_3+W,\qquad
\partial_tW=-U_3+\mathcal{L}W.
\]
Thus, the PV-type component $P$ evolves autonomously, while
$(U_3,W)$ forms a closed wave subsystem. Two representative cases are as follows\begin{enumerate}
    \item When $\mathcal{L}_1=\mathcal{L}$, we can obtain a new set of PV-Wave decomposition unknowns easily\[
\text{PV}_1:=e^{-\mathcal{L} t}\left(\sin\phi\Theta+\cos\phi N\right),\qquad \text{Wave}_1:=e^{-\mathcal{L} t}\left(-\cos\phi\Theta+\sin\phi N\right),\qquad U_{3,1}=e^{-\mathcal{L} t} U_3.
\]
\item If instead $\mathcal{L}_1=0$, then
\[
\partial_tU_3=W,\qquad
\partial_tW=-U_3+\mathcal{L}W,
\]and consequently\[
\partial_t^2U_3-\mathcal{L}\partial_tU_3+U_3=0,\qquad \partial^2_t W-\mathcal{L}\partial_t W+W=0.
\]The operator $\mathcal{L}$ determines the qualitative behaviour of $U_3$ and $W$. Nevertheless, the same change of variables still separates the autonomously evolving PV component from the whole system.
\end{enumerate}

\section{Main results}
The idea of the decomposition method can be naturally extended to study the two systems \eqref{system1} and \eqref{system2}. 
\subsection*{2.1. Boussinesq MHD system}
We start with the first problem \eqref{LBM}. We define the Fourier transform of a function $g$ on $\mathbb{T}^2\times\mathbb{R}$:\[
\widehat{g}(k,\eta,l)=\int_{\mathbb{T}}\int_{\mathbb{T}}\int_{\mathbb{R}} g(x,y,z)e^{-ikx-i\eta y-i lz}\ dxdydz.
\]We also introduce the following mode decomposition in order to categorize different behaviour of each mode.\begin{itemize}
    \item \textbf{(Non-zero modes) ($k\neq 0, \eta\neq0$).}  We define \[
    g_{\neq}(x,y,z):=g(x,y,z)-\frac{1}{|\mathbb{T}|}\int_{\mathbb{T}}g(x,y,z)\ dx-\frac{1}{|\mathbb{T}|}\int_{\mathbb{T}}g(x,y,z)\ dy+\frac{1}{|\mathbb{T}|^2}\int_{\mathbb{T}^2}g(x,y,z)\ dxdy.
    \]
    \item \textbf{(X-zero modes) ($k= 0, \eta\neq0$).} We define\[
    g_{01}(x,y,z):=\frac{1}{|\mathbb{T}|}\int_{\mathbb{T}}g(x,y,z)\ dx-\frac{1}{|\mathbb{T}|^2}\int_{\mathbb{T}^2}g(x,y,z)\ dydx.
    \]
    \item \textbf{(Y-zero modes) ($k\neq 0, \eta=0$).}  We define \[
     g_{10}(x,y,z):=\frac{1}{|\mathbb{T}|}\int_{\mathbb{T}}g(x,y,z)\ dy-\frac{1}{|\mathbb{T}|^2}\int_{\mathbb{T}^2}g(x,y,z)\ dydx.
    \]
    \item \textbf{(Double zero modes) ($k=0, \eta=0$)}.  We define\[
    g_{00}(x,y,z):=\frac{1}{|\mathbb{T}|^2}\int_{\mathbb{T}^2}g(x,y,z)\ dxdy.
    \]
\end{itemize}The Y-zero modes are singled out because, when ($\eta=0$), system \eqref{LBM} decouples into a two-dimensional Boussinesq subsystem and a closed magnetic subsystem, the latter of which may exhibit algebraic growth. As discussed in the introduction, the major difficulty is to analyze the $(u_{3,\neq},\theta_{\neq},b_{3,\neq})$ system. We apply the PV-Wave decomposition method and obtain the following results.

\begin{theorem}(\textbf{Linearized dynamics of \eqref{LBM}}). \label{theo1}Suppose that $\beta\neq0$ and $\alpha>0$, then the following estimates hold
\begin{itemize}
    \item \textbf{Inviscid damping and norm growth of non-zero modes.} The non-zero modes satisfy the following estimate\begin{align*}
         ||(u_3,b_3)_{\neq}||_{L^2}&\lesssim\frac{C_{\alpha,\beta}}{\langle t\rangle}\left(1+\frac{\sqrt{\alpha}}{|\beta|}\right)||(u,\theta,b)_{\neq}(0)||_{H^2}\\
       ||(u_1,u_2,b_1,b_2,\theta)_{\neq}||_{L^2}&\lesssim C_{\alpha,\beta}\left(1+\frac{\sqrt{\alpha}}{|\beta|}\right)\left(1+\frac{1}{|\beta|}\right)||(u,\theta,b)_{\neq}(0)||_{H^1}\\
        ||\theta_{\neq}-\theta_\infty||_{L^2}&\lesssim\frac{\sqrt{\alpha} C_{\alpha,\beta}}{|\beta|\langle t\rangle}\left(1+\frac{\sqrt{\alpha}}{|\beta|}\right)\left(1+\frac{1}{|\beta|}\right)||(u,\theta,b)_{\neq}(0)||_{H^2}.
    \end{align*}Here $C_{\alpha,\beta}:=\left(1+1/\sqrt{\alpha+\beta^2}\right)^\frac{5}{2}e^{\pi/\sqrt{\alpha+\beta^2}}$ and\begin{align*}
    \theta_{\infty}=\theta_{\neq}(0,x-tz,y,z)+\sqrt{\alpha}(\beta\partial_y)^{-1}b_{3,\neq}(0,x-tz,y,z). 
    \end{align*}

    \item \textbf{Inviscid damping, dispersion and instability of Y-zero modes.} The system \eqref{system1} decouples into 2D Boussinesq part and magnetic field part. In particular, there holds\[
    ||(u_2,b_2,b_3)_{10}||_{L^2}+\langle t\rangle^{-1}||b_{1,10}||_{L^2}\lesssim ||(u_2,b)_{10}(0)||_{L^2}.
    \]For $(u_1,u_3,\theta)_{10}$, when $\alpha\neq1/4$, let $r=\text{Re}\sqrt{1/4-\alpha}$, then (\cite{doi:10.1142/S0219891624400058})\[
    \langle t\rangle^{\frac{1}{2}-r}||(u_1,\theta)_{10}||_{L^2\cap L^\infty}+\langle t\rangle^{\frac{3}{2}-r}||u_{3,10}||_{L^2\cap L^\infty}\lesssim||(u_1,u_3,\theta)_{10}(0)||_{H^5}.
    \]When $\alpha=1/4$, then (\cite{doi:10.1142/S0219891624400058})\[
    \langle t\rangle^\frac{1}{2}||(u_1,\theta)_{10}||_{L^2\cap L^\infty}+\langle t\rangle^\frac{3}{2}||u_{3,10}||_{L^2\cap L^\infty}\lesssim\log(2+t)||(u_1,u_3,\theta)_{10}(0)||_{H^5}.
    \]\item \textbf{Dispersive decay of X-zero modes.} $u_{2,01}$ and $u_{3,01}$ experience the following decay\[
    ||(u_2,u_3)_{01}||_{L^\infty}\lesssim\frac{1}{\left(1+\frac{\alpha}{\sqrt{\alpha}+|\beta|}t\right)^\frac{1}{3}}||(u,\theta,b)_{01}(0)||_{W^{3,1}}.
    \]If $\int_{\mathbb{T}}q_{2,l,\text{in}}\ dx=0$, then $\theta_{01},b_{2,01},b_{3,01}$ also experience dispersive decay\[
    ||(\theta,b_2,b_3)_{01}||_{L^\infty}\lesssim\frac{1}{\left(1+\frac{\alpha}{\sqrt{\alpha}+|\beta|}t\right)^\frac{1}{3}}||(u,\theta,b)_{01}(0)||_{W^{3,1}}.
    \]
\end{itemize}    
\end{theorem}\begin{remark} The interaction between the magnetic field and stable stratification generates the norm growth. The sharpness of this growth size can be verified: the linearized PV can be normalized as $\widehat{\Theta}_{\neq}+\frac{\sqrt{\alpha}}{i\beta\eta}\widehat{B}_{3,\neq}=\widehat{\Theta}_{\neq}(0)+\frac{\sqrt{\alpha}}{i\beta\eta}\widehat{B}_{3,\neq}(0)$. Since $B_{3,\neq}\to0$ as $t\to\infty$, we see that $\theta_{\neq}\sim\sqrt{\alpha}/\beta$. When $\sqrt{\alpha}/|\beta|\to 0$, the norm growth of non-zero modes is suppressed, but this will lead to much slower dispersive decay. It is therefore expected that the system is most stable when $\sqrt{\alpha}/|\beta|\sim\mathcal{O}(1)$ at nonlinear level.
\end{remark}
\begin{remark}
     When $y\in\mathbb{R}$, a wave–shear coupled lift-up effect can occur for $u_1,b_1,\theta$ as observed in \cite{yao2026stabilitystratifiedcouette} due to the singular limit $\eta\to0$.
\end{remark}
\begin{remark}
     Our method can be applied to study more general setting where $B^*=\beta(\sigma,1,0)$ for $\sigma\in\mathbb{R}$. The rationality of $\sigma$ can affect the regularity required for inviscid damping, see \cite{liss2020sobolev,wang2026stability} for more details. It also seems possible to apply our method to study the nonlinear stability when $\sqrt{\alpha}/|\beta|\sim 1$.
\end{remark}

\subsection*{2.2. Rotating Boussinesq system}
In the absence of a magnetic field, the rotating Boussinesq system admits a vorticity-dependent potential vorticity \eqref{potentialvorticity}. The normalized linearized PV around the background state \eqref{state2} is given by\[
q_{1,l}=\sqrt{\alpha}\omega_3+(\widetilde{\gamma}-1)\partial_z\theta.
\]The perturbations are defined on $\mathbb{T}\times\mathbb{R}\times\mathbb{T}$ (the case when $z\in\mathbb{R}$ is similar). The corresponding Fourier transform is then defined to be\[
\widehat{g}(k,\eta,l)=\int_{\mathbb{T}}\int_{\mathbb{R}}\int_{\mathbb{T}} g(x,y,z)e^{-ikx-i\eta y-i lz}\ dxdydz
\]Similarly, we introduce the following mode decompositions, \begin{itemize}
    \item \textbf{(Non-zero modes) ($k\neq 0$).} We define \[
    g_{\neq}(x,y,z):=g(x,y,z)-\frac{1}{|\mathbb{T}|}\int_{\mathbb{T}}g(x,y,z)\ dx.
    \]
    \item \textbf{(Simple zero modes) ($k= 0, l\neq0$).} We define\[
    g_{01}(x,y,z):=\frac{1}{|\mathbb{T}|}\int_{\mathbb{T}}g(x,y,z)\ dx-\frac{1}{|\mathbb{T}|^2}\int_{\mathbb{T}^2}g(x,y,z)\ dzdx.
    \]
    \item \textbf{(Double zero modes) ($k=0, l=0$)}. We define\[
    g_{00}(x,y,z):=\frac{1}{|\mathbb{T}|^2}\int_{\mathbb{T}^2}g(x,y,z)\ dxdz.
    \]
\end{itemize}As in the Boussinesq MHD problem, the major difficulty is in analyzing the non-zero mode system $(u_{3,\neq},\theta_{\neq},\omega_{3,\neq})$ system, and the system is analyzed using the PV-Wave decomposition method.

\begin{theorem}(\textbf{Linearized dynamics of \eqref{LRB}}) Suppose that $\alpha>0$, let $\sigma:=\frac{\widetilde{\gamma}(\widetilde{\gamma}-1)}{\alpha}$, we have 
\begin{itemize}
    \item \textbf{Inviscid damping and norm growth of non-zero modes.} The following estimate holds when $\sigma>0$:\begin{align*}
       \langle \sigma\rangle^{-\frac{1}{2}}\langle t\rangle||u_{1,\neq}||_{L^2}+\langle t\rangle||u_{2,\neq}||_{L^2}+||(u_3,\theta)_{\neq}||_{L^2}\lesssim \langle \sigma\rangle^\frac{1}{2}\left(1+\sigma^{-1/2}\right)C_{\alpha,\widetilde{\gamma}}^3e^{C_{\alpha,\widetilde{\gamma}}}||(u,\theta)_{\neq}(0)||_{H^4},
    \end{align*}where $C_{\alpha,\widetilde{\gamma}}=c\left(\frac{1}{\sqrt{\alpha}}+\sqrt{\frac{\widetilde{\gamma}-1}{\widetilde{\gamma}}}+\sqrt{\frac{\widetilde{\gamma}}{\widetilde{\gamma}-1}}+\frac{1}{\sqrt{\widetilde{\gamma}(\widetilde{\gamma}-1)}}\right)$ for some universal constant $c>0$.
     \item \textbf{Spectral instability when $\sigma<0$.} When $\widetilde{\gamma}\in(0,1)$, there exists some initial data and a positive number $\Lambda>0$ such that\[
     ||(u,\theta)_{01}||_{L^2}\gtrsim e^{\Lambda t}.
     \]
\end{itemize}    
\end{theorem}
\begin{remark} The rotating Boussinesq system around the equilibrium \eqref{state2} exhibits strong similarities in the linearized behviour, particularly the inviscid damping pattern, which can be proved in exactly the same way as the case $\widetilde{\gamma}=0$ in \cite{yao2026stabilitystratifiedcouette}. However, the structures of the systems when $\widetilde{\gamma}=0$ and $\widetilde{\gamma}\in(-\infty,0)\cup(1,\infty)$ differ substantially. When $\widetilde{\gamma}=0$, the subsystem \eqref{system2} loses its symmetrizability unless it is reformulated in terms of the linearized potential vorticity, this may incur an additional regularity cost (In \cite{yao2026stabilitystratifiedcouette}, we proved that $H^{\approx4.315}$ regularity was needed for inviscid damping when $\widetilde{\gamma}=0$ and $H^4$ is needed here.) 
\end{remark}

\begin{remark} In contrast to the Boussinesq MHD problem, stable stratification acts primarily as a stabilizing mechanism in the rotating Boussinesq system. The additional growth factor $\sigma^{1/2}$ in $u_{1,\neq}$ is likely an artefact of the method used. To remove this factor may require more precise estimate on $\theta_{\neq}$. Moreover, to prove the nonlinear stability with the inviscid damping of both of the horizontal components of the velocity, it seems that the linearized PV should also be taken into account. However, directly following the method in \cite{yao2026stabilitystratifiedcouette} can fail due to the additional term $\alpha^{-1/2}\widetilde{\gamma}\partial_z\Delta^{-1}_L q_{1,l}$, whose coefficient is no-longer time-integrable. Thus, it seems necessary to consider both $\omega_3$ and $q_{1,l}$ for the nonlinear stability.
\end{remark}
\begin{remark}
    The dispersive estimate for $u_{2,01}$ and $u_{3,01}$ was proved in \cite{huang2026stabilitythreshold3dboussinesq} when $\sigma\neq1$ and $\sigma>0$, we therefore do not present the details here. However, it is worth noting that, similar to the Boussinesq MHD problem, when the simple zero mode of PV has zero initial value, then $u_{1,01}$ and $\theta_{01}$ also experience dispersive decay of the same rate.
\end{remark}
\begin{remark} When $z\in\mathbb{R}$, there is no wave–shear coupled lift-up as observed in \cite{yao2026stabilitystratifiedcouette}.
\end{remark}\vspace{5pt}

This PV-Wave decomposition is particularly effective for systems that admit a compatible symmetrization. When no suitable symmetrizer exists, the method may not be helpful for energy estimate. Section \ref{section5} gives an example. The framework may also extend to other models, such as the full rotating stratified MHD system.\vspace{3pt}

We conclude this section by discussing the combined effects of rotation, stratification, and magnetic fields. As we just saw from the theorems, the interaction of rotation, stratification and magnetic field may produce substantial transient amplification, even though neither mechanism produces the same effect in isolation. This means that the interaction of different types of waves can enrich the dynamics in a non-trivial way. We also record two additional phenomena that are not the main focus of this paper, but are found to be interesting under combined wave effects: 
\begin{itemize}
    \item \textbf{Magnetorotational instability.} This phenomenon has been observed in a variety of settings, see, for instance \cite{sisan2004experimental,hollerbach2005new}. Let us consider the system \eqref{RBM} linearized around the following steady state\[
    v^*=(y,0,0),\qquad \varrho^*=\alpha z,\qquad B^*=(0,0,\beta).
    \]Here $\alpha\geq0$ and $\beta,\widetilde{\gamma}\neq0$. For perturbations posed on $\mathbb{T}\times\mathbb{R}\times\mathbb{T}$, the simple zero mode experiences magnetorotational (spectral) instability when $\widetilde{\gamma}>\beta^2$. The proof can be found in Appendix \ref{MRI}.
    \item \textbf{Enhanced dispersion.} Let us set $\widetilde{\gamma}=0$. Consider the steady state\[
    v^*=(y,0,0),\qquad \varrho^*=\alpha z,\qquad B^*=\beta(\varpi,0,1).
    \]For perturbations set on $\mathbb{R}^2\times\mathbb{T}$, consider the $z$-averaged part. In this case, $(U_3,\Theta,B_3)$ (see Section \ref{notation} and \eqref{notation2} for the definitions) satisfies the following system\[
    \left\{
\begin{aligned}
&\partial_t U_3=\beta\varpi\partial_X B_3+\sqrt{\alpha}\Theta,\\
&\partial_t\Theta=-\sqrt{\alpha}U_3,\\
&\partial_t B_3=\beta\varpi\partial_X U_3.
\end{aligned}
\right.
    \]We can derive a single equation for $U_3$:\[
    \partial^2_t U_3=\beta^2\varpi^2\partial^2_X U_3-\alpha U_3.
    \]This equation is of Klein-Gordon type in the $X$ direction, which produces amplitude decay at rate $\mathcal{O}(t^{-1/2})$ (see \cite{racke1992lectures}, Section 11). In contrast, if either $\alpha$ or $\beta\varpi=0$, there is no amplitude decay\footnote{Recall that 1d wave equation does not experience any decay.}. Thus, the interaction between stratification and the magnetic field enhances dispersion.
\end{itemize}

\subsection*{Notation}\label{notation}
We introduce the different moving frames for the two examples. When $(x,y,z)\in\mathbb{T}^2\times\mathbb{R}$, with the shear flow being $v^*=(z,0,0)$, we write\[
X=x-tz,\qquad Y=y,\qquad Z=z
\]The differential operators in this coordinate system are\begin{align}
       \partial^L_Z:=\partial_Z-t\partial_X,\qquad  \bm{\nabla}_L:=\left(
            \partial_X, \partial_Y, 
            \partial_Z-t\partial_X
        \right)^T,\qquad \Delta_L:=\partial^2_X+\partial^2_Y+(\partial_Z-t\partial_X)^2.\label{notation1}
    \end{align}We write $f$ and $p$ as the Fourier symbols of $-\Delta_H$ and $-\Delta_L$\[
    f:=k^2+\eta^2,\qquad p:=f+(l-kt)^2.
    \]The time derivative is denoted as\[
    \partial_tp:=-2k(l-kt).
    \]When we consider functions defined on $(x,y,z)\in\mathbb{T}\times\mathbb{R}\times\mathbb{T}$, with $v^*=(y,0,0)$, we write\[
    X=x-ty,\qquad Y=y,\qquad Z=z.
    \]The differential operators in this coordinate system are then\begin{align}
       \partial^L_Y:=\partial_Y-t\partial_X,\qquad  \bm{\nabla}_L:=\left(
            \partial_X, \partial_Y-t\partial_X, 
            \partial_Z
        \right)^T,\qquad \Delta_L:=\partial^2_X+(\partial_Y-t\partial_X)^2+\partial^2_Z.\label{notation2}
    \end{align}We write $f$ and $p$ as the Fourier symbols of $-\Delta_H$ and $-\Delta_L$\[
    f:=k^2+(\eta-kt)^2,\qquad p:=f+l^2.
    \]The time derivative is denoted as\[
    \partial_tf=\partial_t p:=-2k(\eta-kt).
    \]
    
    \noindent Define the Japanese bracket for a vector $s\in\mathbb{R}^d$ as $\langle s\rangle:=\sqrt{1+\sum^d_{i=1} s_i^2}$. The Sobolev norm $H^s(\Omega)$ for $\Omega=\mathbb{T}^2\times\mathbb{R}$ and $\Omega=\mathbb{T}\times\mathbb{R}\times\mathbb{T}$ are then defined to be\begin{align*}
&||G||_{H^s(\mathbb{T}^2\times\mathbb{R})}^2=\sum_{(k,\eta)\in\mathbb{Z}^2}\int_{\mathbb{R}} \langle k,\eta,l\rangle^{2s}|\widehat{G}|^2\ dl,\\
&||G||_{H^s(\mathbb{T}\times\mathbb{R}\times\mathbb{T})}^2=\sum_{(k,l)\in\mathbb{Z}^2}\int_{\mathbb{R}} \langle k,\eta,l\rangle^{2s}|\widehat{G}|^2\ d\eta.
\end{align*}

\section{Boussinesq MHD system}

\subsection*{3.1. Non-zero modes}
\begin{figure}[H]
    \centering
    \begin{minipage}[b]{0.49\textwidth}
        \centering
        \includegraphics[width=\textwidth]{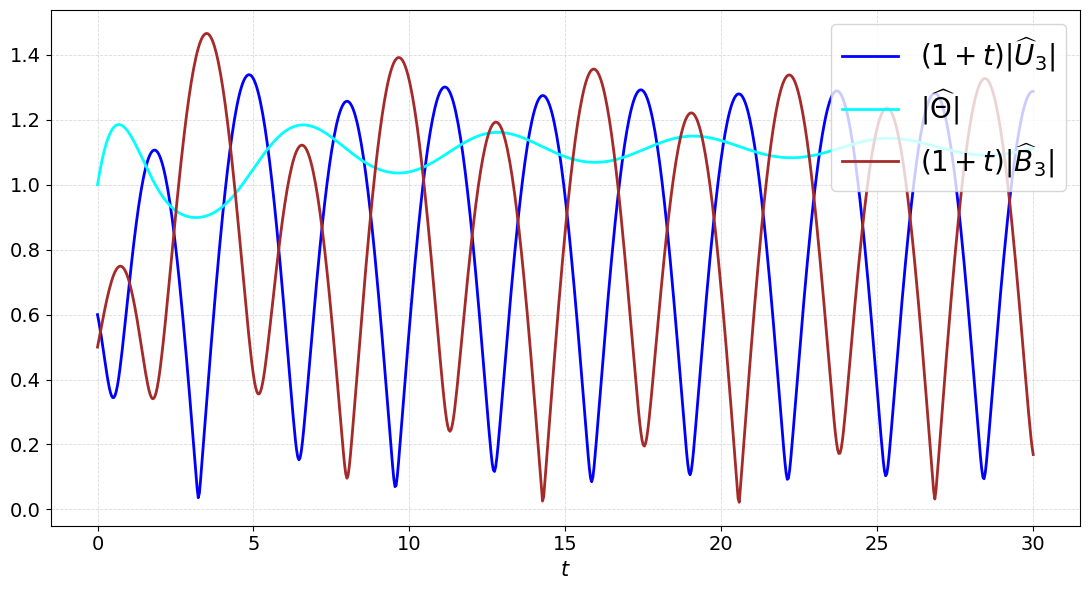}
        
    \end{minipage}
 \hspace{0.01cm}
    \begin{minipage}[b]{0.49\textwidth}
        \centering
        \includegraphics[width=\textwidth]{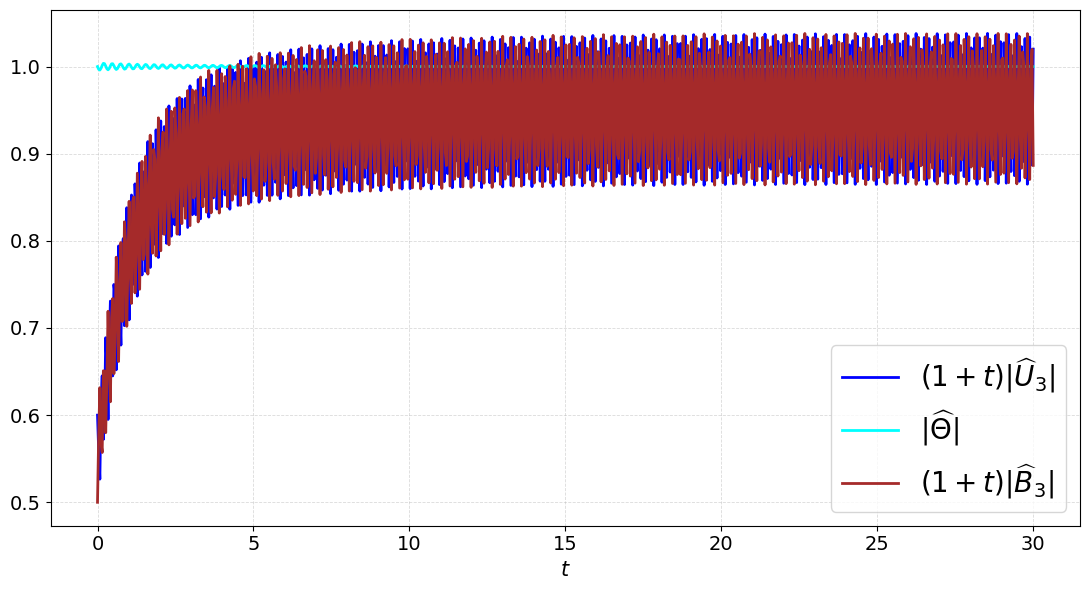}
   
    \end{minipage}
    \hspace{0.01cm}
    \begin{minipage}[b]{0.49\textwidth}
        \centering
        \includegraphics[width=\textwidth]{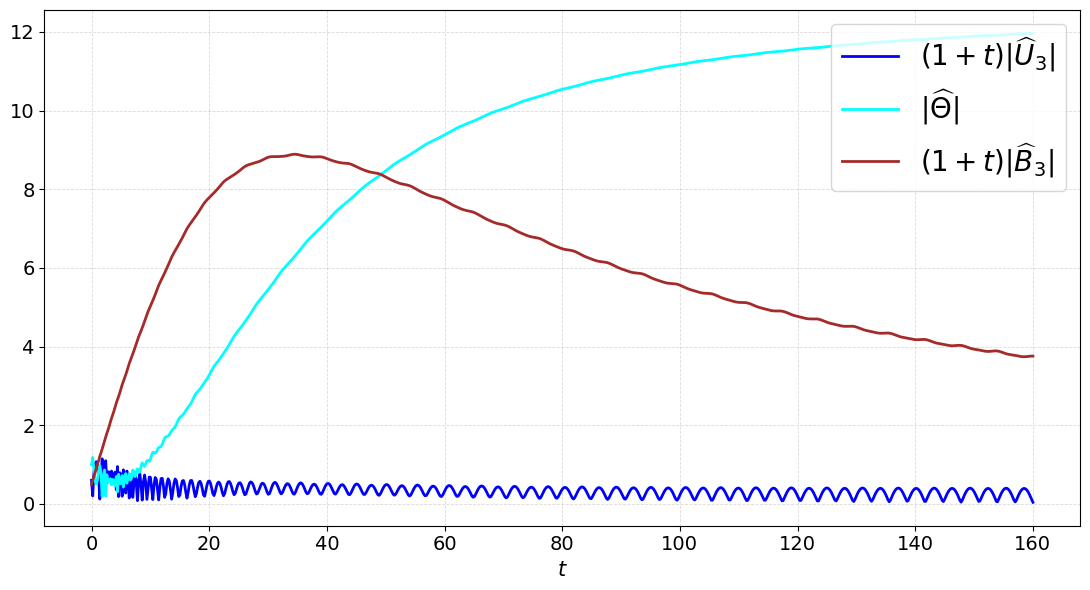}
   
    \end{minipage}
     \hspace{0.01cm}
    \begin{minipage}[b]{0.49\textwidth}
        \centering
        \includegraphics[width=\textwidth]{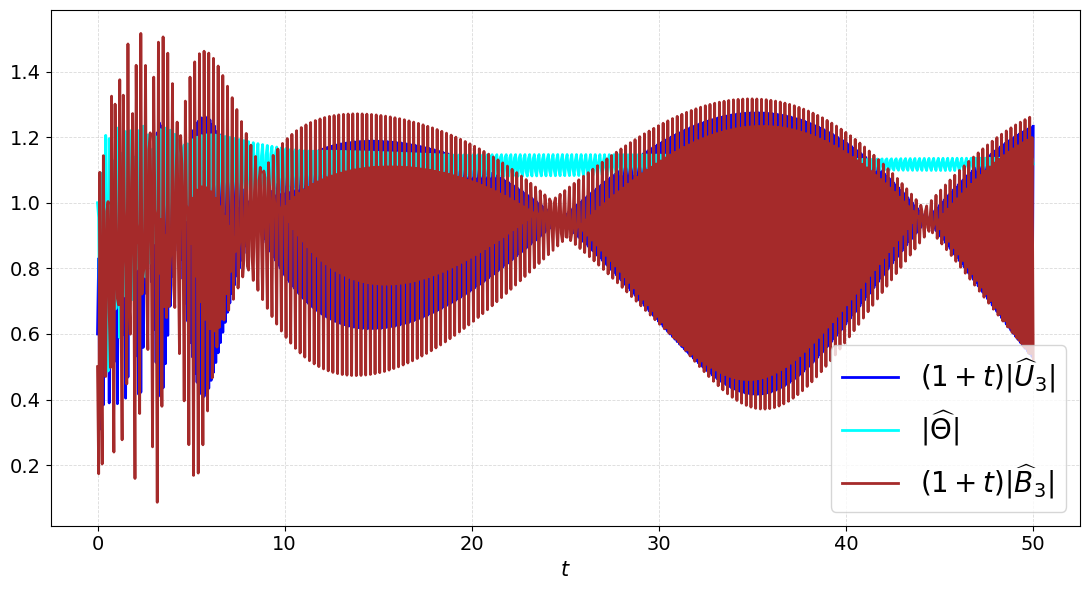}
   
    \end{minipage}
     \caption{Time evolution of $(1+t)|\widehat{U}_3|$, $|\widehat{\Theta}(t)|$, and
$(1+t)|\widehat{B}_3|$ for the fixed Fourier mode
$(k,\eta,l)=(1,1,-1)$.
Top left: $(\alpha,\beta)=(1,1)$, showing bounded rescaled
amplitudes consistent with the decay of $\widehat{U}_3$ and
$\widehat{B}_3$.
Top right: $(\alpha,\beta)=(1,144)$.
Bottom left: $(\alpha,\beta)=(625,1)$, exhibiting amplification
of order $\mathcal{O}(\sqrt{\alpha})$.
Bottom right: $(\alpha,\beta)=(625,25)$, where the magnetic field
suppresses the amplification, consistent with growth of order $\mathcal{O}\!\left(\sqrt{\alpha}/|\beta|\right)$.}\label{fig3}
\end{figure}

In this section, we study the dynamics of non-zero modes. We omit the "$\neq$" subscript for simplicity. Let us recall the linearized system in the moving frame\begin{flalign*}
\left\{
\begin{aligned}
& \partial_t U+U_3\vec{\bm{e}}_1+\bm{\nabla}_L\Delta^{-1}_L(\sqrt{\alpha}\partial_Z^L\Theta-2\partial_X U_3)=\beta\partial_YB+\sqrt{\alpha}\Theta \vec{\bm{e}}_3,\\
& \partial_t\Theta=-\sqrt{\alpha}U_3,\qquad \partial_t B=\beta\partial_YU+B_3\vec{\bm{e}}_1,\\
& \bm{\nabla}_L\cdot U=0,\qquad \bm{\nabla}_L\cdot B=0.
\end{aligned}
\right.
\end{flalign*} $(U_3,\Theta,B_3)$ forms a closed system. After taking the Fourier transform, we obtain the following system for the non-zero modes:\begin{flalign}
\left\{
\begin{aligned}
&\partial_t \widehat{U}_3=-\frac{p'}{p}\widehat{U}_3+i\beta \eta \widehat{B}_3+\sqrt{\alpha}\frac{f}{p}\widehat{\Theta},\\
&\partial_t\widehat{\Theta}=-\sqrt{\alpha}\widehat{U}_3,\\
&\partial_t \widehat{B}_3=i\beta \eta\widehat{U}_3.
\end{aligned}
\right.\label{Xzero}
\end{flalign}Although this system has the same basic structure as \eqref{samplesystem}, its damping behavior is not immediately apparent. We therefore begin with the numerical simulations shown in Figure \ref{fig3}.

\noindent The simulations suggest the following bounds\[
|\widehat{U}_{3}|\sim\frac{1}{t},\qquad  |\widehat{B}_3|\sim\frac{1}{t},\qquad |\widehat{\Theta}|\sim 1.
\]They also suggest norm amplification when $\sqrt{\alpha}/|\beta|\to\infty$. We verify these observations by using the PV-Wave decomposition method.\vspace{9pt}

\noindent\textbf{3.2. The PV-Wave decomposition method.}\vspace{5pt}

\noindent\underline{\textbf{Step 1:} Symmetrizing, identifying the PV unknown and the "wave frequency".}\vspace{5pt}

Consider the following symmetrized variables\[
\widehat{U}_3=\frac{u}{\sqrt{p}},\qquad \widehat{B}_3=\frac{b}{\sqrt{p}},\qquad \widehat{\Theta}=\frac{w}{\sqrt{f}}.
\]The linearized system then becomes\begin{flalign*}
\left\{
\begin{aligned}
&\partial_t u=-\frac{p'}{2p}u+i\beta \eta b+\sqrt{\alpha}\sqrt{\frac{f}{p}}w,\\
&\partial_tw=-\sqrt{\alpha}\sqrt{\frac{f}{p}}u,\\
&\partial_t b=\frac{p'}{2p}b+i\beta \eta u.
\end{aligned}
\right.
\end{flalign*}

To determine the "wave frequency", let us go back to the system \eqref{Xzero}, ignoring the $-\frac{p'}{p}\widehat{U}_3$ (error) term in the $\widehat{U}_3$ equation, taking one more time derivative on $\widehat{U}_3$ equation yields\[
\partial^2_t\widehat{U}_3=-\omega^2\widehat{U}_3+\text{error term}.
\]Here $\omega:=\sqrt{\beta^2\eta^2+\alpha\frac{f}{p}}$. We define the following normalized PV unknown\[
\boxed{Q=w+\frac{\sqrt{\alpha f}}{i\beta\eta\sqrt{p}}b},\qquad \Rightarrow\qquad\partial_t Q=0.
\]\vspace{3pt}

\noindent\underline{\textbf{Step 2:} Determine the corresponding "Wave" unknown.} \vspace{5pt}

 As discussed in the introduction, the aim is to construct an unknown such that $w$ and $b$ can be inverted. We write\[
 Y=r_1 b+r_2w.
 \]The invertibility requires\[
\det\begin{pmatrix}
    r_1&r_2\\
    \frac{\sqrt{\alpha}}{i\beta\eta\sqrt{p}}&\frac{1}{\sqrt{f}}
\end{pmatrix}\neq0,\qquad \iff\qquad i\beta\eta r_1-r_2\sqrt{\alpha}\sqrt{\frac{f}{p}}\neq0.
\]A simple choice is\[
r_1=J(t)i\beta \eta,\qquad r_2=J(t)\sqrt{\alpha}\sqrt{\frac{f}{p}},\qquad J(t)>0
\]where $J$ is to be determined. As a result, we can write\[
\frac{Y}{J}=i\beta \eta b+\sqrt{\alpha}\sqrt{\frac{f}{p}}w.
\]\vspace{3pt}

\noindent\underline{\textbf{Step 3:} Determine the function $J$ for symmetric structure.} \vspace{5pt}

To determine the function $J$, let us first consider equation for $Y$:\begin{align*}
\partial_t Y&=\frac{J'}{J}Y+J\left(-\omega^2 u+\frac{p'}{2p}\left(i\beta\eta b-\sqrt{\alpha}\sqrt{\frac{f}{p}}w\right)\right).
\end{align*}On the other hand, we can invert $b$ and $w$ as\[
b=\frac{i\beta\eta}{\omega^2}\left(\sqrt{\alpha}\sqrt{\frac{f}{p}}Q-\frac{Y}{J}\right),\qquad w=\frac{1}{\omega^2}\left(\sqrt{\alpha}\sqrt{\frac{f}{p}}\frac{Y}{J}+\beta^2\eta^2 Q\right).
\]We therefore write the $u$ equation as\[
\partial_t u=-\frac{p'}{2p}u+\frac{Y}{J}
\]We see that, when $J\omega^2=J^{-1}$, then $(u,Y)$ is symmetric. We therefore take $J=\omega^{-1}$, in this case.\vspace{5pt}

\noindent\underline{\textbf{Step 4:} Constructing an adapted unknown.} \vspace{5pt}

Notice that the term $\frac{p'}{2p}\notin L^1_t$, this can cause some regularity lost.  Motivated by \cite{KNOBEL2026113937}, we consider the following adapted unknown\[
K=u+n(t)Y.
\]This term can be used to cancel the bad terms such as $-\frac{p'}{2p}$ or $\frac{p'}{2p}$. Taking the time derivative, we see that\begin{align*}
    \partial_t(K-nY)+\frac{p'}{2p}(K-nY)&=\omega Y,\\
    \partial_t Y&=-\frac{\omega'}{\omega}Y-\omega(K-nY)\cdots
\end{align*}The first equation becomes\[
\partial_t K-n\left(-\frac{\omega'}{\omega}Y-\textcolor{blue}{\omega K}+\cdots\right)+\textcolor{blue}{\frac{p'}{2p}K}+\cdots=\omega Y
\]Therefore, we see that the bad term can be canceled by setting $n=-\frac{p'}{2p\omega}$. In conclusion, we have defined the following three unknowns\[
\boxed{\underbrace{Q=w+\frac{\sqrt{\alpha f}}{i\beta\eta\sqrt{p}}b}_{\text{PV}},\qquad \underbrace{Y=\frac{i\beta\eta b+\sqrt{\alpha}\sqrt{\frac{f}{p}}w}{\omega}}_{\text{Wave}},\qquad \underbrace{K=u-\frac{p'}{2p\omega} Y}_{\text{Adapted unknown}}.}
\]\vspace{3pt}

\noindent\underline{\textbf{Step 5:} Energy estimate for the new system.}\vspace{5pt}

One can derive the following closed system for $(Q,Y,K)$\begin{flalign*}
\left\{
\begin{aligned}
&\partial_t Y=\frac{\omega'}{\omega}Y-\omega K-\frac{p'\beta^2\eta^2\sqrt{\alpha f}}{p^\frac{3}{2}\omega^3}Q,\\
&\partial_t K=\omega Y-\frac{k^2f}{p^2\omega}Y+\frac{(p')^2\beta^2\eta^2\sqrt{\alpha f}}{2p^\frac{5}{2}\omega^4}Q,\\
&\partial_t Q=0,
\end{aligned}
\right.
\end{flalign*}with the inversion\[
\boxed{b=\frac{i\beta\eta}{\omega^2}\left(\sqrt{\alpha}\sqrt{\frac{f}{p}}Q-\omega Y\right),\qquad w=\frac{1}{\omega^2}\left(\omega\sqrt{\alpha}\sqrt{\frac{f}{p}}Y+\beta^2\eta^2 Q\right)}.
\]Now we can derive the energy estimate, consider the energy functional\[
\mathbb{E}:=\frac{1}{2}\left[|Y|^2+|K|^2\right].
\]We have\[
\partial_t\mathbb{E}\leq 2\cdot\bm{1}_{\omega'\geq0}\frac{\omega'}{\omega}\mathbb{E}+\frac{k^2f}{p^2\omega}\mathbb{E}+2\sqrt{2}\left(1+\frac{1}{\sqrt{\alpha+\beta^2}}\right)\frac{|p'|\beta^2\eta^2\sqrt{\alpha f}}{p^\frac{3}{2}\omega^3}\sqrt{\mathbb{E}}|Q|.
\]Here we have used the fact that $\frac{|p'|}{p\omega}\leq\frac{2}{\sqrt{\alpha+\beta^2}}$. Since\[
\int_{\mathbb{R}}\frac{k^2f}{p^2\omega}\ ds\leq \int_{\mathbb{R}} \frac{k^2}{\sqrt{\alpha+\beta^2}p}\ ds\leq\frac{\pi}{\sqrt{\alpha+\beta^2}},
\]direct integration yields\[
\sqrt{\mathbb{E}}\lesssim e^{\int_0^t\bm{1}_{\omega'\geq0}\frac{\omega'}{\omega}\ ds}\left(\sqrt{\mathbb{E}(0)}+\left(1+\frac{1}{\sqrt{\alpha+\beta^2}}\int_0^t\frac{|p'|\beta^2\eta^2\sqrt{\alpha f}}{(p\omega^2)^\frac{3}{2}}\right)|Q|\right).
\]Due to the conservation of $Q$, we have \[
|Q|\lesssim\left(1+\frac{\sqrt{\alpha f}}{|\beta\eta|\sqrt{p(0)}}\right)\left(|w(0)|+|b(0)|\right).
\]Moreover, the following integral is finite, \[
I:=\int_0^t\frac{|p'|\beta^2\eta^2\sqrt{\alpha f}}{p^\frac{3}{2}\omega^3}\ ds\leq\int_{\mathbb{R}}\frac{2|k(l-kt)|\beta^2\eta^2\sqrt{\alpha f}}{(\beta^2\eta^2 p+\alpha f)^\frac{3}{2}}\ ds.
\]Applying the change of variable $l-ks=\sqrt{k^2+\eta^2}\tau=\sqrt{f}\tau$, we see that\begin{equation}
I\leq2\int\frac{\beta^2\eta^2\sqrt{\alpha}f^\frac{3}{2}|\tau|}{(\beta^2\eta^2 f(1+\tau^2)+\alpha f)^\frac{3}{2}}\ d\tau\leq2\int\frac{\beta^2\eta^2\sqrt{\alpha}|\tau|}{(\beta^2\eta^2(1+\tau^2)+\alpha)^\frac{3}{2}}\ d\tau.\label{integralestimate1}
\end{equation}Let $\beta\eta\tau=\sqrt{\alpha+\beta^2\eta^2}t$, we have\begin{equation}
I\leq\frac{2\sqrt{\alpha}}{\sqrt{\alpha+\beta^2\eta^2}}\int\frac{|t|}{(1+t^2)^\frac{3}{2}}\ dt\lesssim1.\label{boundmmm}
\end{equation}Finally, consider the following product\begin{align*}
e^{\int \bm{1}_{\omega'\geq0}\frac{\omega'}{\omega}\ ds}\left(1+\frac{\sqrt{\alpha f}}{|\beta\eta|\sqrt{p}(0)}\right)&\leq\sup_m\frac{\omega(m)}{\omega(0)}\left(1+\frac{\sqrt{\alpha f}}{|\beta\eta|\sqrt{p}(0)}\right)\\&=\sup_m\sqrt{\frac{\beta^2\eta^2p(m)+\alpha f}{p(m)}}\cdot\sqrt{\frac{p(0)}{\beta^2\eta^2p(0)+\alpha f}}\left(1+\frac{\sqrt{\alpha f}}{|\beta\eta|\sqrt{p(0)}}\right)\\
&\lesssim1+\frac{\sqrt{\alpha}}{|\beta\eta|}.
\end{align*}We can thus conclude that\[
\mathbb{E}\lesssim  \left(1+\frac{\alpha}{\beta^2}\right) \left(1+\frac{1}{\sqrt{\alpha+\beta^2}}\right)e^{\frac{2\pi}{\sqrt{\alpha+\beta^2}}}\left(\mathbb{E}(0)+|b(0)|^2+|w(0)|^2\right).
\]
\noindent\underline{\textbf{Step 6:} Energy estimate for original unknowns.}\vspace{5pt}

Observe the relations:\begin{align*}
  |Y|+|K|&\lesssim\left(1+\frac{1}{\sqrt{\alpha+\beta^2}}\right)\left(|b|+|w|+|u|\right),\\
  |b|+|w|+|u|&\lesssim \left(1+\frac{1}{\sqrt{\alpha+\beta^2}}\right) \left( |Y|+|K|\right)+ \left(1+\frac{\sqrt{\alpha}}{|\beta|}\right)|w(0)|+|b(0)|.
\end{align*}We therefore obtain the bound\[
\sqrt{p}|\widehat{U}_3|+\sqrt{p}|\widehat{B}_3|+\sqrt{f}|\widehat{\Theta}|\lesssim\left(1+\frac{\sqrt{\alpha}}{|\beta|}\right) C_{\alpha,\beta}\left(\sqrt{p(0)}|\widehat{U}_3(0)|+\sqrt{p}(0)|\widehat{B}(0)|+\sqrt{f}|\widehat{\Theta}(0)|\right).
\]Here $C_{\alpha,\beta}:=\left(1+\frac{1}{\sqrt{\alpha+\beta^2}}\right)^\frac{5}{2}e^{\frac{\pi}{\sqrt{\alpha+\beta^2}}}$. Hence, we obtain the inviscid damping estimate\[
\langle t\rangle||u_{3,\neq}||_{L^2}+\langle t\rangle||b_{3,\neq}||_{L^2}+||\theta_{\neq}||_{L^2}\lesssim \left(1+\frac{\sqrt{\alpha}}{|\beta|}\right) C_{\alpha,\beta} ||(u_3,\theta,b_3)_{\neq}(0)||_{H^2}.
\]\hfill \qedsymbol
\begin{remark}
    Direct estimates of $u$, $b$, and $w$ do not yield the desired decay estimate, not only because of the presence of the non-integrable (in time) term $\frac{p'}{2p}$, but also because the norm growth cannot be precisely quantified. A further reason for using the PV-Wave decomposition method is that direct symmetrization yields the wrong scaling and therefore fails to recover the correct inviscid damping estimate for $(u_3,\theta,b_3)$.
\end{remark}

\begin{remark}
We can invert the linearized PV $Q$ to deduce the long-time behaviour of $\theta$ since $b_{3}\to0$ as $t\to\infty$:\[
    \theta\to  \theta_{\infty}:=\theta(0,x-tz,y,z)+\sqrt{\alpha}(\beta\partial_y)^{-1}b_3(0,x-tz,y,z)\qquad\text{as }t\to\infty. 
    \]
\end{remark}
\vspace{8pt}

\noindent\textbf{3.3. Analysis of $(u_1,u_2,b_1,b_2)_{\neq}$.} Again, motivated by \cite{KNOBEL2026113937}, we consider the following adapted unknown\[
\widetilde{V}=U+(\beta\partial_Y)^{-1}B_3\vec{\bm{e}}_1.
\]In this case, the linearized system turns to\begin{flalign}
\left\{
\begin{aligned}
& \partial_t \widetilde{V}+\bm{\nabla}_L\Delta^{-1}_L(\sqrt{\alpha}\partial_Z^L\Theta-2\partial_X\widetilde{V}_3)=\beta\partial_YB+\sqrt{\alpha}\Theta\vec{\bm{e}}_3,\\
& \partial_t\Theta=-\sqrt{\alpha}\widetilde{V}_3,\\
& \partial_t B=\beta\partial_Y\widetilde{V},\\
& \bm{\nabla}\cdot \widetilde{V}=\partial_X(\beta\partial_Y)^{-1}B_3,\qquad \bm{\nabla}\cdot B=0.
\end{aligned}
\right.
\label{LRSM2}
\end{flalign}We consider the modified vorticity, i.e.\[
\widetilde{W}_3:=\partial_X \widetilde{V}_2-\partial_Y\widetilde{V}_1
\]and the current density $J_3$, these two unknowns form a simple system\begin{flalign}
\left\{
\begin{aligned}
& \partial_t \widetilde{W}_3=\beta \partial_Y J_3,\\
&\partial_t J_3=\beta\partial_Y \widetilde{W}_3.
\end{aligned}
\right.
\label{LRSM3}
\end{flalign}The system can be solved explicitly. Alternatively, by standard energy estimate, we have\[
\left|\widehat{\widetilde{W}_3}\right|^2+\left|\widehat{J}_3\right|^2=\left|\widehat{\widetilde{W}_3}(0)\right|^2+\left|\widehat{J}_3(0)\right|^2
\]For the velocity field, we have\begin{flalign*}
\left\{
\begin{aligned}
&\partial_x u_2-\partial_y u_1=\widetilde{\omega}_3+\frac{b_3}{\beta},\\
&\partial_x u_1+\partial_y u_2=-\partial_z u_3.
\end{aligned}
\right.
\label{LRSM3}
\end{flalign*}Thus, we obtain the following estimate\[
||(u_1,u_2,b_1,b_2)_{\neq}||_{L^2}\lesssim\left(1+\frac{\sqrt{\alpha}}{|\beta|}\right)\left(1+\frac{1}{|\beta|}\right)C_{\alpha,\beta}||(u,\theta,b)_{\neq}(0)||_{H^1}.
\]
\subsection*{3.4. Y-zero modes}
In the $Y$-zero-mode system, the transport effect induced by the shear flow persists. However, the main difference from the non-zero mode system is that the magnetic field is decoupled from the velocity and temperature fields due to resonance. The linearized system is then given by (the subscript "10" is omitted here.)\begin{flalign*}
\left\{
\begin{aligned}
& \partial_t U+U_3\bm{\vec{e}}_1+\bm{\nabla}_{X,Z}\Delta^{-1}_L(\sqrt{\alpha}\partial_Z^L\Theta-2\partial_X U_3)=\sqrt{\alpha}\Theta \bm{\vec{e}}_3,\\
& \partial_t\Theta=-\sqrt{\alpha}U_3,\qquad \partial_t B=B_3\bm{\vec{e}}_1,\\
& \partial_X U_1+\partial_Z^L U_3=0,\qquad \partial_X B_1+\partial_{Z}^L B_3=0.
\end{aligned}
\right.
\end{flalign*}$U_2$ and $B$ can be computed explicitly:\[
U_2=U_2(0),\qquad B_1=B_1(0)+t B_3(0),\qquad B_2=B_2(0),\qquad B_3=B_3(0).
\]Notice that $(U_1,U_3,\Theta)$ reduces to the 2D Boussinesq system, whose solution can be solved explicitly at the linear level in terms of the Whittaker functions \cite{doi:10.1142/S0219891624400058}.

\subsection*{3.5. X-zero modes}
In the $x$-averaged part, the transport effect due to the shear flow disappears. Omitting the subscript "01", the Fourier-transformed system is formulated as\begin{equation}
\left\{
\begin{aligned}
&\partial_t\widehat{u}_{1}+\widehat{u}_3=i\beta\eta\widehat{b}_1,\qquad \partial_t \widehat{u}_2=i\beta\eta\widehat{b}_2-\sqrt{\alpha}\frac{\eta l}{\eta^2+l^2}\widehat{\theta},\\
&\partial_t\widehat{u}_3=i\beta\eta\widehat{b}_3+\sqrt{\alpha}\frac{\eta^2}{\eta^2+l^2}\widehat{\theta},\qquad \partial_t\widehat\theta=-\sqrt{\alpha}\widehat{u}_3,\\
&\partial_t\widehat{b}_1=i\beta\eta\widehat{u}_1+\widehat{b}_3,\qquad \partial_t\widehat{b}_2=i\beta \eta\widehat{u}_2,\qquad \partial_t\widehat{b}_3=i\beta\eta\widehat{u}_3.
\end{aligned}
\right.
\end{equation}The above ODE system can be solved explicitly. Define the following phase functions\[
\Phi_1=\beta\eta,
\qquad
\Phi_2=
\sqrt{\beta^2\eta^2+\frac{\alpha\eta^2}{\eta^2+l^2}},
\]the solutions are then computed as{\footnotesize\[
\begin{aligned}
\widehat u_1(t,\eta,l)
={}&
\widehat u_1(0,\eta,l)\cos(\Phi_1 t)
+i\widehat b_1(0,\eta,l)\sin(\Phi_1 t)
-\frac{\widehat u_3(0,\eta,l)}{\Phi_2}\sin(\Phi_2 t)
\\
&-\frac{\sqrt{\alpha}\eta^2}
 {(\eta^2+l^2)\Phi_2^2}
 \widehat\theta(0,\eta,l)
 \left[1-\cos(\Phi_2 t)\right]-\frac{i\widehat b_3(0,\eta,l)}{\Phi_1}
\left[
\cos(\Phi_1 t)-1
+\frac{\Phi_1^2}{\Phi_2^2}
 \left(1-\cos(\Phi_2 t)\right)
\right],
\\[1ex]
\widehat u_2(t,\eta,l)
={}&
\widehat u_2(0,\eta,l)\cos(\Phi_2 t)
+\frac{
 i\Phi_1\widehat b_2(0,\eta,l)
 -\displaystyle\sqrt{\alpha}\frac{\eta l}{\eta^2+l^2}
  \widehat\theta(0,\eta,l)}
 {\Phi_2}
 \sin(\Phi_2 t),
\\[1ex]
\widehat u_3(t,\eta,l)
={}&
\widehat u_3(0,\eta,l)\cos(\Phi_2 t)
+\frac{
 i\Phi_1\widehat b_3(0,\eta,l)
 +\displaystyle\sqrt{\alpha}\frac{\eta^2}{\eta^2+l^2}
  \widehat\theta(0,\eta,l)}
 {\Phi_2}
 \sin(\Phi_2 t),
\\[1ex]
\widehat\theta(t,\eta,l)
={}&
\widehat\theta(0,\eta,l)
-\frac{\sqrt{\alpha}}{\Phi_2}
 \widehat u_3(0,\eta,l)\sin(\Phi_2 t)
-\frac{i\sqrt{\alpha}\Phi_1}{\Phi_2^2}
 \widehat b_3(0,\eta,l)
 \left[1-\cos(\Phi_2 t)\right]-\frac{\alpha\eta^2}
 {(\eta^2+l^2)\Phi_2^2}
 \widehat\theta(0,\eta,l)
 \left[1-\cos(\Phi_2 t)\right],
\\[1ex]
\widehat b_1(t,\eta,l)
={}&
\widehat b_1(0,\eta,l)\cos(\Phi_1 t)
+i\widehat u_1(0,\eta,l)\sin(\Phi_1 t)
+\frac{\widehat b_3(0,\eta,l)}{\Phi_1}\sin(\Phi_1 t),
\\[1ex]
\widehat b_2(t,\eta,l)
={}&
\widehat b_2(0,\eta,l)
+\frac{i\Phi_1}{\Phi_2}
 \widehat u_2(0,\eta,l)\sin(\Phi_2 t)
-\frac{\Phi_1^2}{\Phi_2^2}
 \widehat b_2(0,\eta,l)
 \left[1-\cos(\Phi_2 t)\right]-\frac{i\sqrt{\alpha}\eta l\,\Phi_1}
 {(\eta^2+l^2)\Phi_2^2}
 \widehat\theta(0,\eta,l)
 \left[1-\cos(\Phi_2 t)\right],
\\[1ex]
\widehat b_3(t,\eta,l)
={}&
\widehat b_3(0,\eta,l)
+\frac{i\Phi_1}{\Phi_2}
 \widehat u_3(0,\eta,l)\sin(\Phi_2 t)
-\frac{\Phi_1^2}{\Phi_2^2}
 \widehat b_3(0,\eta,l)
 \left[1-\cos(\Phi_2 t)\right]+\frac{i\sqrt{\alpha}\eta^2\Phi_1}
 {(\eta^2+l^2)\Phi_2^2}
 \widehat\theta(0,\eta,l)
 \left[1-\cos(\Phi_2 t)\right].
\end{aligned}
\]}\begin{remark}
    When $y\in\mathbb{T}$, $\eta$ is bounded away from $0$ (except the $\eta=0$ mode) and hence the classical lift-up effect arising in the Navier Stokes setting (see, e.g. \cite{9e0bb666-1e4e-3aa2-87b0-f916c2336491}) is suppressed. However, when $y\in\mathbb{R}$, low frequency limit $\eta\to0$ can lead to the wave-shear coupled lift-up effect, as observed in \cite{yao2026stabilitystratifiedcouette}. In particular, $\lim_{\eta\to0}\frac{\sin(\Phi_i t)}{\Phi_i}=t$ for fixed $l$, and we see that $u_1,b_1,\theta$ can experience this type of lift-up.
\end{remark}\vspace{3pt}

Since $y\in\mathbb{T}$, the Fourier variable $\eta$ is discrete. Moreover, $\Phi_1=\beta\eta$ is independent of the continuous Fourier variable $l$, so the corresponding branch does not exhibit dispersive decay. By contrast, for each fixed non-zero $\eta\in\mathbb{Z}$, the phase $\Phi_2$ depends nontrivially on $l\in\mathbb{R}$, and this dependence gives rise to dispersive decay. The explicit formulas show that $u_2$ and $u_3$ are governed entirely by this dispersive branch and hence exhibit dispersive decay. From the explicit solution, we see that $u_2$ and $u_3$ are purely dispersive and experience dispersive decay. Moreover, if the initial X-zero modes of PV is zero, then $\theta_{01},b_{2,01},b_{3,01}$ are also dispersive. To study the decay rate, we follow the framework developed in \cite{coti2024stability} and provide a simplified proof. \begin{lemma}
 Let $\mathcal{L}:=\partial_y\sqrt{\beta^2+\alpha|\bm{\nabla}_{y,z}|^{-2}}$, the following estimate holds\[
 \left\|e^{t\mathcal{L}}f\right\|_{L^\infty\left(\mathbb{T}\times\mathbb{R}\right)}\lesssim\frac{1}{\left(1+\frac{\alpha}{\sqrt{\alpha}+|\beta|}t\right)^\frac{1}{3}}||f||_{W^{3,1}}.
 \]  
\end{lemma}\begin{proof}Let\begin{gather*}
f(y,z)=\frac{1}{2\pi}\sum_{\eta\neq0}\widehat{f}_\eta(\eta,z)e^{i\eta y},\qquad \hat{f}_\eta(\eta,z)=\int_{\mathbb{T}}f(y,z)e^{-i\eta y}\ dy,\\\widehat{f}(\eta,l)=\int_{\mathbb{R}} \widehat{f}_\eta(\eta,z)e^{-ilz}\ dz.
\end{gather*}We define the following Littlewood-Paley decomposition\begin{flalign*}
\phi\in C^\infty(\mathbb{R}),\qquad \phi(l)=\left\{
\begin{aligned}
& 1&|l|\leq 1\\
& 0&|l|>\frac{4}{3}
\end{aligned}
\right.,\qquad \varphi(l)=\phi\left(\frac{l}{2}\right)-\phi(l),\qquad \varphi_j(l)=\varphi(2^{-j}l).
\end{flalign*}Denote the projection operator as $\mathbb{P}_j\psi:=\mathbb{F}^{-1}\left(\varphi_j(l)\widehat{\psi}\right)$, we therefore have\begin{align}
\left\|e^{t\mathcal{L}}\widehat{f}_\eta\right\|_{L^\infty_z}&\leq\sum_{j\in\mathbb{Z}}\left\|\mathcal{F}^{-1}_z\left(e^{it\eta\sqrt{\beta^2+\alpha(\eta^2+l^2)^{-1}}}\varphi_j\varphi_j^1\widehat{f}(\eta,l)\right)\right\|_{L^\infty_z}\nonumber\\
    &\lesssim\sum_{j\in\mathbb{Z}}\left\|\mathcal{F}_z^{-1}\left(e^{it\eta\sqrt{\beta^2+\alpha(\eta^2+l^2)^{-1}}}\varphi_j\right)\right\|_{L^\infty}\left\|\mathbb{P}^1_j\widehat{f}_\eta(\eta,\cdot)\right\|_{L^1_z}.\label{asodyasodasoudsd}
\end{align}Here $\varphi^1_j$ and $\varphi_j$ are chosen to have similar supports such that $\varphi_j=\varphi_j\varphi_j^1$. It remains to study the following integral after writing $l=\eta\zeta$\[
I=\eta\int_{\mathbb{R}} e^{it\Phi(\zeta)}\varphi(2^{-j}\eta\zeta)\ d\zeta.
\]The phase function $\Phi$ is defined as follows\[
\Phi(\zeta)=\sqrt{\beta^2\eta^2+\frac{\alpha}{1+\zeta^2}}+\Omega\zeta=\sqrt{\alpha}\left(\sqrt{\lambda+\frac{1}{1+\zeta^2}}+\widetilde{\Omega}\zeta\right)=:\sqrt{\alpha}\phi(\zeta),\qquad \Omega=\frac{\eta z}{t}
\]Here $\lambda:=\beta^2\eta^2/\alpha>0$ and $\widetilde{\Omega}=\alpha^{-1/2}\Omega$. To obtain a bound for $I$, we compute the derivatives of $\phi$:\[
\begin{aligned}
\phi'(\zeta)
&=
\frac{
\widetilde{\Omega}(1+\zeta^2)^{3/2}
\sqrt{1+\lambda(1+\zeta^2)}-\zeta
}{
(1+\zeta^2)^{3/2}
\sqrt{1+\lambda(1+\zeta^2)}
},
\\[2mm]
\phi''(\zeta)
&=
\frac{
\lambda(1+\zeta^2)(3\zeta^2-1)+2\zeta^2-1
}{
(1+\zeta^2)^{5/2}
\left[1+\lambda(1+\zeta^2)\right]^{3/2}
}.
\end{aligned}
\]Solving the equation $\phi''=0$, we see that there can be two degenerate critical points $\zeta^2_\pm=\frac{\sqrt{(1+\lambda)(1+4\lambda)}-(1+\lambda)}{3\lambda}\in[1/3,1/2]$. We consider for simplicity the case where $\eta\geq 0$ and split the integral $I$ into\[
I=I_++I_-,\qquad I_+=\eta\int_{\mathbb{R}^+} e^{it\Phi(\zeta)}\varphi(2^{-j}\eta\zeta)\ d\zeta.
\]We focus on the $I_+$ in the rest of the proof. Notice that\[
\zeta_+^*\in [2^{j-2}\eta^{-1},2^{j+2}\eta^{-1}]\iff j\in[j_0,j_0+4],\qquad j_0:=\log_2\eta-\frac{5}{2}.
\]We consider the following two cases:\begin{itemize}
    \item When $j\in[j_0,j_0+4]$. The critical point can be degenerate, we therefore split the integral further:\[
    I_+^1+I_+^2:=\eta\int_{B_\epsilon(\zeta_+)}+\eta\int_{B^c_\epsilon(\zeta_+)}e^{it\Phi(\zeta)}\varphi(2^{-j}\eta\zeta)\ d\zeta
    \]for some constant $\epsilon>0$ to be determined. Notice that $|\phi''|\gtrsim \frac{1+|\lambda|}{1+|\lambda|^\frac{3}{2}}\epsilon$. For the second integral, we can apply the van der Corput lemma, and deduce that\[
|I_+^2|\lesssim|\eta|\left(\frac{1+|\lambda|}{1+|\lambda|^\frac{3}{2}}\epsilon\sqrt{\alpha}t\right)^{-\frac{1}{2}}.
    \]Together with the simple estimate $|I_+^1|\lesssim|\eta|\epsilon$ by picking $\epsilon^
\frac{3}{2}=\left(\frac{1+|\lambda|}{1+|\lambda|^\frac{3}{2}}\sqrt{\alpha}t\right)^{-\frac{1}{2}}$, we obtain the bound\[
|I_+|\lesssim\frac{|\eta|^\frac{4}{3}}{\left(1+\frac{\alpha}{\sqrt{\alpha}+|\beta|}t\right)^\frac{1}{3}}.
\]\item When $j\notin[j_0,j_0+4]$, if $j<j_0$, we have both $|I_+|\lesssim 2^j$ and $|I_+|\lesssim|\eta|\left(\frac{1+|\lambda|}{1+|\lambda|^\frac{3}{2}}\sqrt{\alpha}t\right)^{-1/2}$ due to the van der Corput lemma. If $j>j_0+4$, we have $|\zeta|\lesssim 2^j|\eta|^{-1}$, hence $|\phi''|\gtrsim |\zeta|^{-4}\frac{1+|\lambda|}{1+|\lambda|^\frac{3}{2}}\gtrsim 2^{-4j}\eta^{4}\frac{1+|\lambda|}{1+|\lambda|^\frac{3}{2}}$. Again, applying van der Corput lemma, we obtain the bound $|I_+|\lesssim 2^{2j}|\eta|^{-1}\left(\frac{1+|\lambda|}{1+|\lambda|^\frac{3}{2}}\sqrt{\alpha}t\right)^{-1/2}$.
\end{itemize}Summing over $j\in\mathbb{Z}$ yields \[
\left\|e^{t\mathcal{L}}\widehat{f}_\eta\right\|_{L^\infty_z}\lesssim\frac{1}{\left(1+\frac{\alpha}{\sqrt{\alpha}+|\beta|}t\right)^\frac{1}{3}}\left(|\eta|^{-1/2}||\widehat{f}_\eta||_{W^{2+\varepsilon,1}}+|\eta|^{4/3}||\widehat{f}_\eta||_{W^{\varepsilon,1}}\right).
\]Summing over $\eta$, and for small $\varepsilon\in(0,1/5)$, we have\[
\left\|e^{t\mathcal{L}}f\right\|_{L^\infty}\lesssim \frac{1}{\left(1+\frac{\alpha}{\sqrt{\alpha}+|\beta|}t\right)^\frac{1}{3}}||f||_{W^{1/2+\varepsilon,1}_yW_z^{2+\varepsilon,1}\cap W^{7/3+\varepsilon,1}_y\cap W^{\varepsilon,1}_z}\lesssim \frac{1}{\left(1+\frac{\alpha}{\sqrt{\alpha}+|\beta|}t\right)^\frac{1}{3}}||f||_{W^{3,1}}.
\]\end{proof}
\begin{remark}
   The estimate above shows that the decay becomes slower in the regime $|\beta|\gg\sqrt{\alpha}$. This can be seen from the phase function $\Phi_2$: as $\beta\to\infty$, $\Phi_2\to|\beta\eta|$, which does not produce strong oscillation in the phase space. This again indicates that the interaction of the two dispersive mechanisms does not necessarily enhance decay.
\end{remark}

\subsection*{3.6. Double-zero modes}
The double-zero modes evolve more simply. In this setting, the two divergence-free conditions reduce to\[
\partial_z u_{3,00}=0,\qquad \partial_z b_{3,00}=0.
\]Therefore $u_{3,00}$ and $b_{3,00}$ can depend only on $t$. However, if $u_{3,00}(0)$ and $b_{3,00}(0)$ are assumed to be in $L^2_z$, they must be identically equal to zero. Thus, the linearized system for double-zero modes becomes\begin{flalign*}
\left\{
\begin{aligned}
& \partial_t u_{00}+\partial_zp_{00}\vec{\bm{e}}_3=\sqrt{\alpha}\theta_{00}\vec{\bm{e}}_3,\\
& \partial_t\theta_{00}=0,\qquad  \partial_t b_{00}=0.
\end{aligned}
\right.
\end{flalign*}We can immediately deduce that all components are time-independent, and we have\[
||(u,\theta,b)_{00}||_{L^2}\lesssim||(u,\theta,b)_{00}(0)||_{L^2}.
\]

\section{Rotating Boussinesq system}
    In the absence of a magnetic field, the linearized system around the steady state \eqref{state2}, with the subscript "$\neq$" omitted is given by\begin{flalign}
\left\{
\begin{aligned}
&\partial_t u+y\partial_x u+u_2\vec{\bm{e}}_1+\bm{\nabla}\Delta^{-1}\left(-2\partial_x u_2+\widetilde{\gamma}\omega_3+\sqrt{\alpha}\partial_z\theta\right)+\widetilde{\gamma}\vec{\bm{e}}_3\times u=\sqrt{\alpha}\theta\vec{\bm{e}}_3,\\
&\bm{\nabla}\cdot u=0,\\
    &\partial_t\theta+y\partial_x\theta=-\sqrt{\alpha} u_3.
\end{aligned}
\right.\label{linearRB}
\end{flalign}The setting here shares some similarities with the non-rotating Boussinesq problem studied in \cite{yao2026stabilitystratifiedcouette}. The difference is that $(\omega_3,u_3,\theta)$ now forms a symmetrizable system. In the moving frame, the system is formulated as\begin{flalign}
\left\{
\begin{aligned}
&\partial_t  \widehat{U}_3=\frac{f'l^2}{pf}\widehat{U}_3+\frac{i\widetilde{\gamma}l}{p}\widehat{W}_3+\sqrt{\alpha}\frac{f}{p}\widehat{\Theta}-\frac{2ik^2l}{pf}\widehat{W}_3,\\
&\partial_t\widehat{W}_3=-il(1-\widetilde{\gamma})\widehat{U}_3,\\
    &\partial_t\widehat{\Theta}=-\sqrt{\alpha}\widehat{U}_3.
\end{aligned}
\right.\label{2342342}
\end{flalign}Similarly, we use numerical simulations to identify the expected decay behavior and examine whether large parameter values lead to norm growth. Since $W_3$ can be treated as an auxiliary variable, it suffices to analyze the behaviour of $(U_3,\Theta)$. To derive a closed system for these two unknowns, we introduce the linearized PV\[
S:=\sqrt{\alpha}W_3+(\widetilde{\gamma}-1)\partial_z\Theta,\qquad \Rightarrow\qquad \partial_t S=0.
\]By writing $\sigma:=\widetilde{\gamma}(\widetilde{\gamma}-1)/\alpha$, we have\begin{flalign}
\left\{
\begin{aligned}
&\partial_t\Theta=-\sqrt{\alpha} U_3,\\
 &\partial_t U_3=2\partial_{XZ}\Delta^{-1}_L\Delta_H^{-1}\left[\frac{\partial_X S+(1-\widetilde{\gamma})\partial_{XZ}\Theta}{\sqrt{\alpha}}-\partial_{YZ}^LU_3\right]+\sqrt{\alpha}\Delta_H\Delta^{-1}_L\theta+\sqrt{\alpha}\sigma\partial^2_Z\Delta^{-1}_L\Theta-\frac{\widetilde{\gamma}}{\sqrt{\alpha}}\Delta^{-1}_L\partial_Z S.
\end{aligned}
\right.\label{mkmkmkmk}
\end{flalign}\begin{figure}[H]
    \centering
    \begin{minipage}[b]{0.49\textwidth}
        \centering
        \includegraphics[width=\textwidth]{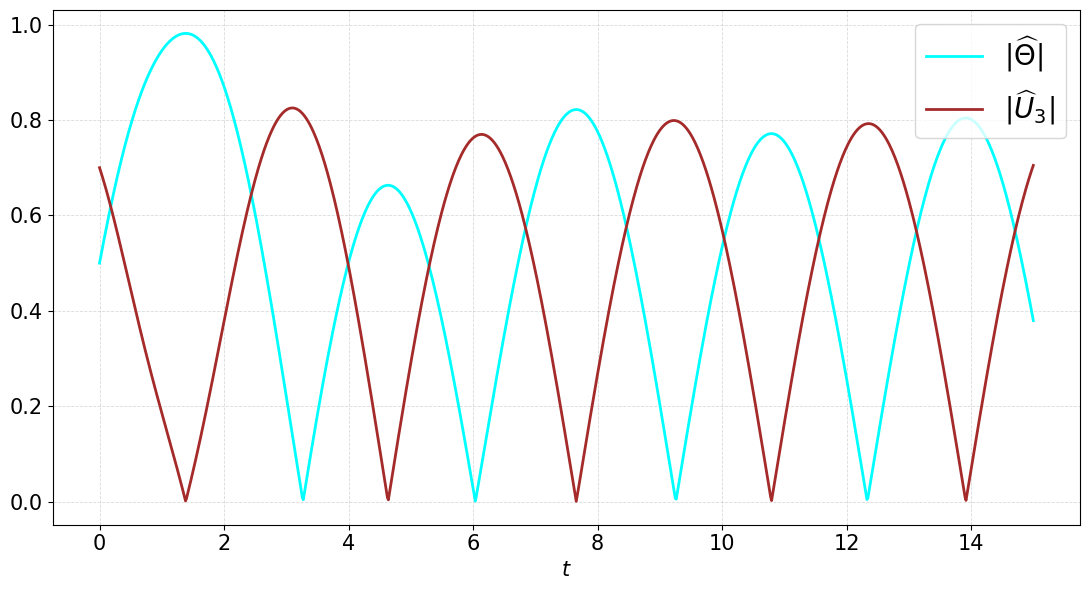}
        
    \end{minipage}
    \hspace{0.01cm}
    \begin{minipage}[b]{0.49\textwidth}
        \centering
        \includegraphics[width=\textwidth]{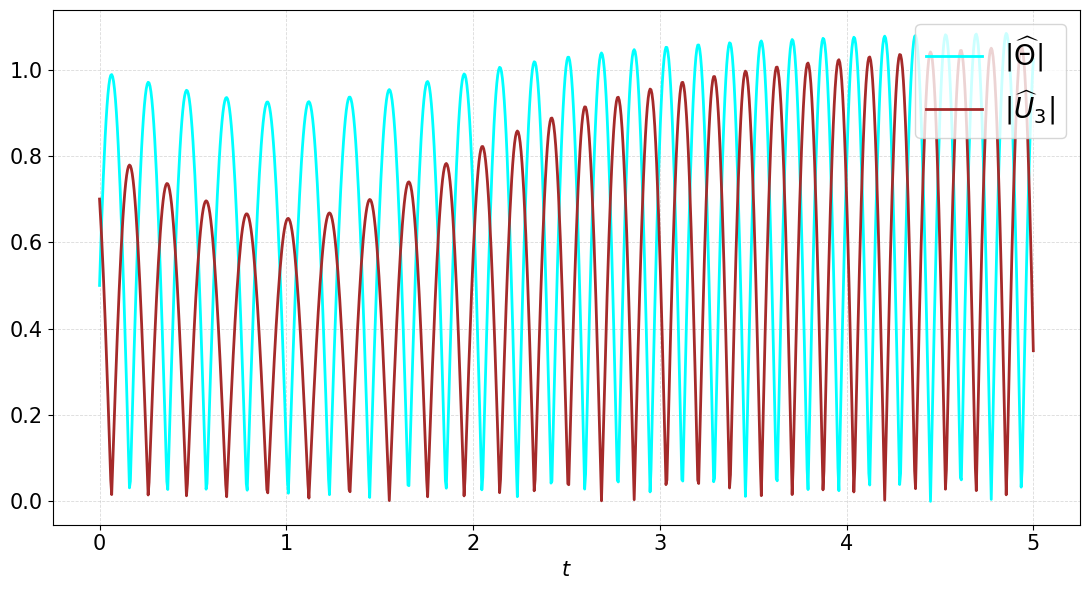}
   
    \end{minipage}
    \hspace{0.01cm}
    \begin{minipage}[b]{0.49\textwidth}
        \centering
        \includegraphics[width=\textwidth]{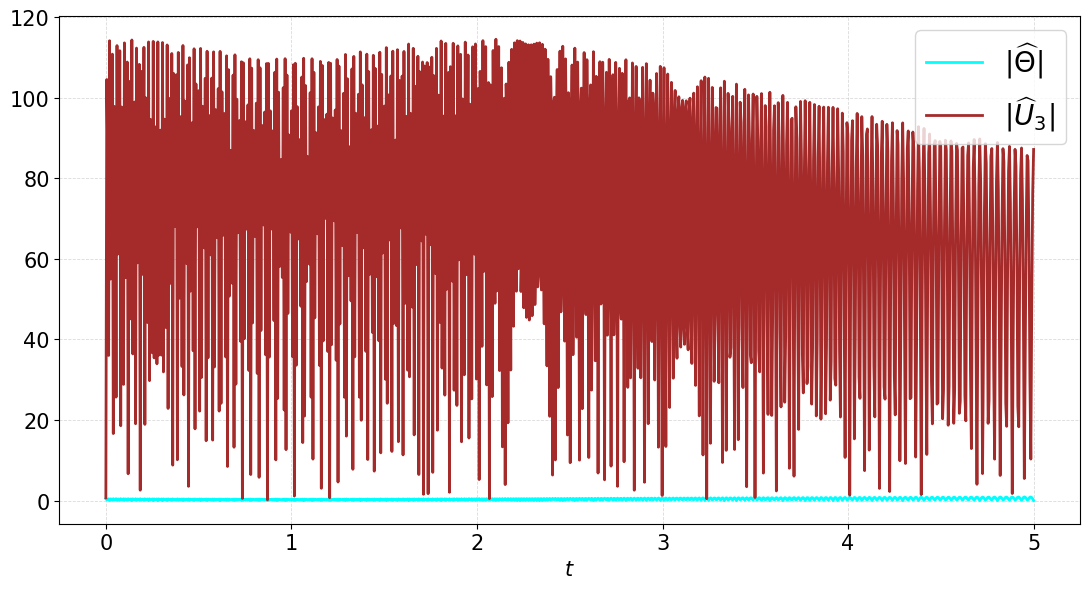}
   
    \end{minipage}
    \hspace{0.01cm}
    \begin{minipage}[b]{0.49\textwidth}
        \centering
        \includegraphics[width=\textwidth]{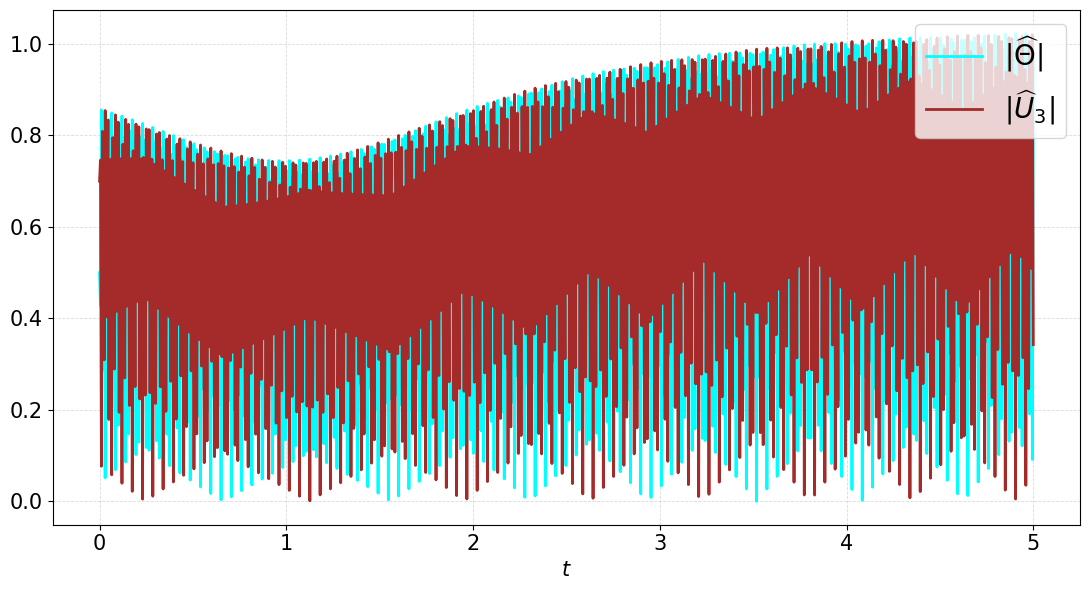}
   
    \end{minipage}
     \caption{
Time evolution of $\left|\widehat{\Theta}\right|$ and $\left|\widehat{U}_3\right|$
for the fixed Fourier mode $(k,\eta,l)=(1,1,-1)$. Top left: $(\alpha,\widetilde{\gamma})=(1,2)$.
Top right: $(\alpha,\widetilde{\gamma})=(400,2)$.
Bottom left: $(\alpha,\widetilde{\gamma})=(1,400)$.
Bottom right:
$(\alpha,\widetilde{\gamma})=(160000,400)$.
}\label{fig2}
\end{figure}

\noindent The numerical simulations shown in Figure \ref{fig2} suggest behavior similar to that in the non-rotating case $\widetilde{\gamma}=0$ studied in \cite{yao2026stabilitystratifiedcouette}:\[
|U_3|\sim 1,\qquad |\Theta|\sim 1.
\]Moreover, it turns out that the solution exhibits norm growth when rotation is relatively large compared to the stratification strength.\vspace{5pt}

We now verify these observations. Notice that the dynamics are extremely simple when $l=0$:\[
\partial_t\Theta=-\sqrt{\alpha} U_3.\qquad \partial_t U_3=\sqrt{\alpha}\Theta.
\]Hence we consider the case where $l\neq0$, and apply the PV-Wave decomposition method.\vspace{9pt}

\noindent\textbf{4.1. The PV-Wave decomposition method.}\vspace{5pt}

\noindent\underline{\textbf{Step 1:} Symmetrizing, identifying the PV unknown and the "wave frequency".}\vspace{5pt}

We first determine the wave frequency, in \eqref{2342342}, we ignore the "error terms", i.e. $\frac{f'l^2}{pf}\widehat{U}_3$ and $-\frac{2ik^2l}{pf}\widehat{W}_3$, we therefore have\[
\partial^2_t\widehat{U}_3=-\omega^2\widehat{U}_3+\text{error term},\qquad\omega^2=\frac{\alpha f+\widetilde{\gamma}(\widetilde{\gamma}-1)l^2}{p}.
\]We set\[
\widehat{\Theta}=A,\qquad \widehat{U}_3=\underbrace{\sqrt{\frac{f}{p}}}_{=:b}B,\qquad \widehat{W}_3=\underbrace{\sqrt{\frac{\widetilde{\gamma}-1}{\widetilde{\gamma}}}\sqrt{f}}_{=:c}C,
\]the linearized system then becomes\begin{flalign*}
\left\{
\begin{aligned}
&\partial_t  A=-\sqrt{\alpha}bB,\\
&\partial_tB=\frac{f'l^2}{pf}B+\frac{i\widetilde{\gamma}l}{p}\frac{c}{b}C-\frac{b'}{b}B+\sqrt{\alpha}\frac{f}{p}\frac{1}{b}A-\frac{2ik^2l}{pf}\frac{c}{b}C,\\
    &\partial_tC=-il(1-\widetilde{\gamma})\frac{b}{c}B-\frac{c'}{c}C.
\end{aligned}
\right.
\end{flalign*}We define the PV mode (assume $\widetilde{\gamma}(\widetilde{\gamma}-1)>0$)\[
K_R=\sqrt{\frac{\widetilde{\gamma}-1}{\widetilde{\gamma}}}\sqrt{f}C+\frac{\widetilde{\gamma}-1}{\sqrt{\alpha}}ilA=\sqrt{\frac{\widetilde{\gamma}-1}{\widetilde{\gamma}}}\left(\sqrt{f}C+\text{sgn}(\widetilde{\gamma})i\sqrt{\sigma}l A\right),\qquad \partial_t K_R=0.
\]\vspace{3pt}

\noindent\underline{\textbf{Step 2:} Determine the corresponding "Wave" unknown.} \vspace{5pt}

Set the "Wave" unknown to be\[
Y:=r_1C+r_2A.
\]Again, to ensure that $A$ and $C$ can be inverted via $K_R$ and $Y$, we require\[
\det\begin{pmatrix}
    r_1&r_2\\\sqrt{\frac{\widetilde{\gamma}-1}{\widetilde{\gamma}}}\sqrt{f}&\frac{\widetilde{\gamma}-1}{\sqrt{\alpha}}il
\end{pmatrix}\neq0.
\]We choose\[
r_1=\frac{\widetilde{\gamma}-1}{\sqrt{\alpha}}il J(t),\qquad r_2=\sqrt{\frac{\widetilde{\gamma}-1}{\widetilde{\gamma}}}\sqrt{f} J(t).
\]Thus, we have\[
Y=J(t)\left(\frac{\widetilde{\gamma}-1}{\sqrt{\alpha}}il C+\sqrt{\frac{\widetilde{\gamma}-1}{\widetilde{\gamma}}}\sqrt{f}A\right).
\]Thus, $A$ and $C$ can be inverted:\begin{align*}
    A=\frac{\alpha\widetilde{\gamma}}{\widetilde{\gamma}-1}\frac{\sqrt{\frac{\widetilde{\gamma}-1}{\widetilde{\gamma}}}\sqrt{f}\frac{Y}{J}-\frac{\widetilde{\gamma}-1}{\sqrt{\alpha}}il K_R}{p\omega^2},\qquad C=-\frac{\alpha\widetilde{\gamma}}{\widetilde{\gamma}-1}\frac{\frac{\widetilde{\gamma}-1}{\sqrt{\alpha}}il\frac{Y}{J}-\sqrt{\frac{\widetilde{\gamma}-1}{\widetilde{\gamma}}}\sqrt{f}K_R}{p\omega^2},
\end{align*}\vspace{3pt}

\noindent\underline{\textbf{Step 3:} Determine the function $J$ for symmetric structure.} \vspace{5pt}

Direct computation yields\[
\partial_t Y=\frac{J'}{J}Y+J\left(-\frac{1}{\sqrt{\alpha}}\sqrt{\frac{\widetilde{\gamma}-1}{\widetilde{\gamma}}}\sqrt{p}\omega^2 B+\frac{f'}{2f}\left(\sqrt{\frac{\widetilde{\gamma}-1}{\widetilde{\gamma}}}\sqrt{f}A-\frac{\widetilde{\gamma}-1}{\sqrt{\alpha}}il C\right)\right).
\]Similarly, we have\[
\partial_t B=\frac{f'l^2}{2pf}B+\frac{\sqrt{\alpha}}{\sqrt{p}}\sqrt{\frac{\widetilde{\gamma}}{\widetilde{\gamma}-1}}\frac{Y}{J}-\frac{2ik^2l}{\sqrt{p}f}\sqrt{\frac{\widetilde{\gamma}-1}{\widetilde{\gamma}}}C.
\]Hence, we need to pick $J(t)=\frac{\sqrt{\alpha}}{\omega\sqrt{p}}\sqrt{\frac{\widetilde{\gamma}}{\widetilde{\gamma}-1}}$ to ensure the symmetric structure of the system. We can then derive the following closed system \begin{flalign}
\left\{
\begin{aligned}
&\partial_t  B=\frac{f'l^2}{2fp}B+\omega Y-\frac{2(\widetilde{\gamma}-1)k^2l^2}{fp\omega}Y-\frac{2i\alpha k^2l}{\sqrt{pf}p\omega^2} K_R,\\
&\partial_tY=-\omega B-\frac{\widetilde{\gamma}(\widetilde{\gamma}-1)l^2f'}{2fp\omega^2}Y-\frac{i\alpha\widetilde{\gamma}lf'}{\sqrt{f}(p\omega^2)^\frac{3}{2}}K_R,\\
    &\partial_tK_R=0.
\end{aligned}
\right.\label{pvwave1}
\end{flalign}\vspace{5pt}

\noindent\underline{\textbf{Step 4:} Constructing an adapted unknown.} \vspace{5pt}

To cancel the $\frac{f'l^2}{2fp}B$ term in \eqref{pvwave1}, we set the following ansatz\[
Z=B+n(t)Y,
\]To deduce $n$, we observe that\[
\partial_t Y=-\omega Z+\cdots.
\]On the other hand, we can deduce an equation for $Z$ from $B$,\begin{align*}
&\partial_t Z-n'Y-n\partial_t Y=\frac{f'l^2}{2pf}(Z-nY)+\cdots\\
\Rightarrow\qquad&\partial_t Z-n' Y-n(\textcolor{blue}{-\omega Z}+\cdots)=\frac{f'l^2}{2pf}\left(\textcolor{blue}{Z}+\cdots\right).
\end{align*}As a result, cancellation can be done by picking $n(t)=\frac{f'l^2}{2pf\omega}$. Consequently, $(Z,K_R,Y)$ satisfy the following system\begin{flalign}
\left\{
\begin{aligned}
&\partial_t  Z=\omega Y-\frac{2(\widetilde{\gamma}-1)k^2l^2}{fp\omega}Y+\frac{\partial_t\left(\frac{f'l^2}{2pf}\right)-\left(\frac{f'l^2}{2pf}\right)^2}{\omega}Y-\left(\frac{2i\alpha k^2l}{\sqrt{pf}p\omega^2}+\frac{i\alpha\widetilde{\gamma}lf'}{\omega(p\omega^2)^\frac{3}{2}\sqrt{f}}\frac{f'l^2}{2pf}\right)K_R,\\
&\partial_tY=-\omega Z+\frac{\omega'}{\omega}Y-\frac{i\alpha\widetilde{\gamma}lf'}{\sqrt{f}(p\omega^2)^\frac{3}{2}}K_R,\\
    &\partial_tK_R=0.
\end{aligned}
\right.\label{pvwave2}
\end{flalign}\vspace{3pt}

\noindent\underline{\textbf{Step 5:} Energy estimate for the new system.}\vspace{5pt}  

We consider the energy functional\[
\mathbb{E}:=\frac{1}{2}\left[|Z|^2+|Y|^2\right],
\]we have\begin{align*}
\partial_t\mathbb{E}&\leq 2\cdot\bm{1}_{\omega'\geq0}\frac{\omega'}{\omega}\mathbb{E}+\frac{\sqrt{2}\alpha|\widetilde{\gamma}||l||f'|}{\sqrt{f}(p\omega^2)^\frac{3}{2}}|K_R|\sqrt{\mathbb{E}}+\frac{2|\widetilde{\gamma}-1|k^2l^2}{fp\omega}\mathbb{E}\\&+\frac{1}{\omega}\left|\partial_t\left(\frac{f'l^2}{2pf}\right)-\left(\frac{f'l^2}{2pf}\right)^2\right|\mathbb{E}+\left(\frac{4\alpha k^2|l|}{\sqrt{pf}p\omega^2}+\frac{\alpha|\widetilde{\gamma}lf'|}{\omega(p\omega^2)^\frac{3}{2}\sqrt{f}}\frac{|f'|l^2}{2pf}\right)|K_R|\sqrt{\mathbb{E}}.
\end{align*}Notice that\begin{align*}
\int_{\mathbb{R}}\frac{|\widetilde{\gamma}-1|k^2l^2}{fp\omega}\leq\sqrt{\frac{\widetilde{\gamma}-1}{\widetilde{\gamma}}}\int_{\mathbb{R}}\frac{k^2}{f}&\leq\pi\sqrt{\frac{\widetilde{\gamma}-1}{\widetilde{\gamma}}},\\\int_{\mathbb{R}}\frac{1}{\omega}\left|\partial_t\left(\frac{f'l^2}{2pf}\right)\right|+\int_{\mathbb{R}}\frac{1}{\omega}\left(\frac{f'l^2}{2pf}\right)^2&\lesssim\frac{1}{\min\{\sqrt{\alpha},\sqrt{\widetilde{\gamma}(\widetilde{\gamma}-1)}\}}.
\end{align*}For the first term, we have\[
\text{exp}\left(\int_0^t\bm{1}_{\omega'\geq 0}\frac{\omega'}{\omega}\right)\leq\sup_{t_1,t_2\in\mathbb{R}}\sqrt{\frac{f(t_1)+\sigma l^2}{p(t_1)}\cdot\frac{p(t_2)}{f(t_2)+\sigma  l^2}}\lesssim\sqrt{\sigma+\sigma^{-1}}.
\]To deal with the last term, notice that\begin{align*}
|K_R|\sqrt{\mathbb{E}}&\leq\sqrt{\frac{\widetilde{\gamma}-1}{\widetilde{\gamma}}}\sqrt{f(0)}|C(0)|\sqrt{\mathbb{E}}+\sqrt{\sigma}\sqrt{\frac{\widetilde{\gamma}-1}{\widetilde{\gamma}}}|l||A(0)|\sqrt{\mathbb{E}}\\
&\lesssim\sqrt{\frac{\widetilde{\gamma}-1}{\widetilde{\gamma}}}\left(f(0)|C(0)|^2+\sqrt{\sigma}l^2|A(0)|^2+(1+\sqrt{\sigma})\mathbb{E}\right).
\end{align*}We therefore only need to consider the following three integrals\begin{align*}
    I_1&=\int\frac{\alpha k^2|l|}{\sqrt{pf}p\omega^2}\sqrt{\frac{\widetilde{\gamma}-1}{\widetilde{\gamma}}}(1+\sqrt{\sigma}),\qquad I_2=\int\frac{\alpha|\widetilde{\gamma}lf'|}{\omega(p\omega^2)^\frac{3}{2}\sqrt{f}}\frac{|f'|l^2}{pf}\sqrt{\frac{\widetilde{\gamma}-1}{\widetilde{\gamma}}}(1+\sqrt{\sigma})\\
    I_3&=\int\frac{\alpha|\widetilde{\gamma}lf'|}{\sqrt{f}(p\omega^2)^\frac{3}{2}}\sqrt{\frac{\widetilde{\gamma}-1}{\widetilde{\gamma}}}(1+\sqrt{\sigma}).
\end{align*}We have\begin{align*}
    I_1&\leq\int\frac{|k|}{f}\sqrt{\frac{\widetilde{\gamma}-1}{\widetilde{\gamma}}}\lesssim\sqrt{\frac{\widetilde{\gamma}-1}{\widetilde{\gamma}}}\\
    I_2&\lesssim\frac{\sqrt{\sigma}(1+\sqrt{\sigma})}{\sqrt{\alpha}}\int\frac{k^2l^2}{(f+\sigma l^2)^2}\lesssim\left(\frac{1}{\sqrt{\alpha}}+\frac{1}{\sqrt{\widetilde{\gamma}(\widetilde{\gamma}-1)}}\right)\int\frac{k^2}{f}\lesssim \frac{1}{\sqrt{\alpha}}+\frac{1}{\sqrt{\widetilde{\gamma}(\widetilde{\gamma}-1)}}\\
    I_3&\lesssim\int\frac{\sqrt{\sigma}(1+\sqrt{\sigma})|kl|}{(f+\sigma l^2)^\frac{3}{2}}\lesssim\frac{1+\sqrt{\sigma}}{\sqrt{1+\sigma l^2}}\lesssim 1.
\end{align*}We can therefore deduce the following estimate\[
\mathbb{E}\lesssim (\sigma+\sigma^{-1}) C_{\alpha,\widetilde{\gamma}}e^{C_{\alpha,\widetilde{\gamma}}}\left(\mathbb{E}(0)+p(0)\left(|A(0)|^2+|C(0)|^2\right)\right),
\]where $C_{\alpha,\widetilde{\gamma}}=c\left(\frac{1}{\sqrt{\alpha}}+\sqrt{\frac{\widetilde{\gamma}-1}{\widetilde{\gamma}}}+\sqrt{\frac{\widetilde{\gamma}}{\widetilde{\gamma}-1}}+\frac{1}{\sqrt{\widetilde{\gamma}(\widetilde{\gamma}-1)}}\right)$ for some universal constant $c>0$. \vspace{5pt}

\noindent\underline{\textbf{Step 6:} Energy estimate for original unknowns.}\vspace{5pt}

We can deduce the following relations\begin{align*}
|A|+|B|+|C|&\lesssim C_{\alpha,\widetilde{\gamma}}\left(p^\frac{1}{2}(0)\left(|A(0)|+|C(0)|\right)+|Y|+|Z|\right),\\
|Y|+|Z|&\lesssim C_{\alpha,\widetilde{\gamma}}\left(|A|+|B|+|C|\right)
\end{align*}As a result, we obtain the following estimate\[
|A|+|B|+|C|\lesssim \langle\sigma\rangle^\frac{1}{2}\left(1+\sigma^{-1/2}\right)C^\frac{5}{2}_{\alpha,\widetilde{\gamma}}e^{C_{\alpha,\widetilde{\gamma}}}p^\frac{1}{2}(0)\left(|A(0)|+|B(0)|+|C(0)|\right).
\]From above, we deduce the estimate\[
||\theta_{\neq}||_{L^2}+||u_{3,\neq}||_{L^2}\lesssim \langle\sigma\rangle^\frac{1}{2}\left(1+\sigma^{-1/2}\right)C^\frac{5}{2}_{\alpha,\widetilde{\gamma}}e^{C_{\alpha,\widetilde{\gamma}}}\left\|(u,\theta)_{\neq}(0)\right\|_{H^2}.
\]To deduce the inviscid damping of $u_{1,\neq}$ and $u_{2,\neq}$, we use the fact that\begin{align*}
   & \partial_XU_1+\partial_Y^L U_2=-\partial_Z U_3\\
    &\partial_X U_2-\partial_Y^L U_1= W_3.
\end{align*}For $U_2$, we obtain that\[
|U_2|\lesssim\frac{|k|}{f}|W_3|+\frac{|l|}{f^\frac{1}{2}}|U_3|,
\]from which we deduce that\[
||u_{2,\neq}||_{L^2}\lesssim \frac{\langle\sigma\rangle^\frac{1}{2}}{\langle t\rangle}\left(1+\sigma^{-1/2}\right)C^3_{\alpha,\widetilde{\gamma}}e^{C_{\alpha,\widetilde{\gamma}}}\left\|(u,\theta)_{\neq}(0)\right\|_{H^4}.
\]For $U_1$, observe that\[
|W_3|\lesssim |W_3(0)|+\sqrt{\frac{\widetilde{\gamma}-1}{\widetilde{\gamma}}}\sqrt{\sigma}|l|\left(|A(0)|+|A|\right),
\]we therefore obtain a bit rougher bound\[
||u_{1,\neq}||_{L^2}\lesssim \frac{\langle\sigma\rangle}{\langle t\rangle}\left(1+\sigma^{-1/2}\right)C^3_{\alpha,\widetilde{\gamma}}e^{C_{\alpha,\widetilde{\gamma}}}\left\|(u,\theta)_{\neq}(0)\right\|_{H^4}.
\]\hfill \qedsymbol

\begin{remark}
    Indeed, the scaling for $W_3$ in the symmetrization is incorrect, since $u_1\sim t^{-1}$, the correct scaling should be $W_3\sim 1$ instead of $W_3\sim t$. 
\end{remark}
\begin{remark}
    The boundedness and inviscid damping can be proved in a much simpler way, to be more specific, by considering the following energy functional for system \eqref{mkmkmkmk}:\[
\mathbb{E}:=\frac{1}{2}\left(|\Theta|^2+\frac{p}{f+\sigma l^2}|U_3|^2\right).
\]However, the drawback is that the norm growth size cannot be deduced as precise as the above method.
\end{remark}
\subsection*{4.2. Simple zero modes}

The system for the simple zero modes simplfies to
\[\left\{
\begin{aligned}
&\partial_tu_{1,01}=(\widetilde{\gamma}-1)u_{2,01},\\
    &\partial_tu_{2,01}=-\sqrt{\alpha}\partial_{yz}\Delta^{-1}_{y,z}\theta_{01}-\widetilde{\gamma}\partial^2_z\Delta^{-1}u_{1,01},\\
    &\partial_t u_{3,01}=\sqrt{\alpha}\partial^2_y\Delta^{-1}\theta_{01}+\widetilde{\gamma}\partial_{yz}\Delta^{-1}u_{1,01},\\
    &\partial_t\theta_{01}=-\sqrt{\alpha} u_{3,01}.
\end{aligned}
\right.\]Again, the solution can be computed explicitly in Fourier variables. We can obtain a single equation for $\widehat{u}_{3,01}$\[
\partial^2_t\widehat{u}_3=-\frac{\alpha\eta^2+\widetilde{\gamma}(\widetilde{\gamma}-1)l^2}{\eta^2+l^2}\widehat{u}_3.
\]In general, let $\widetilde{\gamma}\in(-\infty,0)\cup(1,\infty)$ and denote the phase function by $\Phi_1:=\sqrt{\frac{\alpha\eta^2+\widetilde{\gamma}(\widetilde{\gamma}-1)l^2}{\eta^2+l^2}}$, the system is solved as follows
\begin{equation}
\left\{
\begin{aligned}\widehat{u}_{1,01}=&\frac{\alpha\eta^2\widehat{u}_{1,01}(0)-\sqrt{\alpha}(\widetilde{\gamma}-1)\eta l\widehat{\theta}_{01}(0)}{\alpha\eta^2+\widetilde{\gamma}(\widetilde{\gamma}-1)l^2}\\&+(\widetilde{\gamma}-1)\left[\widehat{u}_{2,01}(0)\frac{\sin(\Phi_1t)}{\Phi_1}+\frac{\sqrt{\alpha}\eta l\widehat{\theta}_{01}(0)+\widetilde{\gamma}l^2\widehat{u}_{1,01}(0)}{\alpha\eta^2+\widetilde{\gamma}(\widetilde{\gamma}-1)l^2}\cos(\Phi_1 t)\right],\\
    \widehat{u}_{2,01}=&\widehat{u}_{2,01}(0)\cos(\Phi_1 t)-\frac{\eta l\sqrt{\alpha}\widehat{\theta}_{01}(0)+\widehat{\gamma}l^2\widehat{u}_{1,01}(0)}{\sqrt{\eta^2+l^2}\sqrt{\alpha\eta^2+\widetilde{\gamma}(\widetilde{\gamma}-1)l^2}}\sin(\Phi_1 t),\\
\widehat{u}_{3,01}=&\widehat{u}_{3,01}(0)\cos(\Phi_1 t)+\frac{\sqrt{\alpha}\eta^2\widehat{\theta}_{01}(0)+\widetilde{\gamma}\eta l\widehat{u}_{1,01}(0)}{\sqrt{\eta^2+l^2}\sqrt{\alpha\eta^2+\widetilde{\gamma}(\widetilde{\gamma}-1)l^2}}\sin(\Phi_1 t),\\
\widehat{\theta}_{01}=&\frac{\widetilde{\gamma}(\widetilde{\gamma}-1)l^2\widehat{\theta}_{01}(0)-\sqrt{\alpha}\widetilde{\gamma}\eta l\widehat{u}_{1,01}(0)}{\alpha\eta^2+\widetilde{\gamma}(\widetilde{\gamma}-1)l^2}-\frac{\sqrt{\alpha}}{\Phi_1}\widehat{u}_{3,01}(0)\sin(\Phi_1 t)+\frac{\alpha\eta^2\widehat{\theta}_{01}(0)+\sqrt{\alpha}\widetilde{\gamma}\eta l\widehat{u}_{1,01}(0)}{\alpha\eta^2+\widetilde{\gamma}(\widetilde{\gamma}-1)l^2}\cos(\Phi_1 t).
\end{aligned}
\right.\label{explicitsolution000}
\end{equation}
Both the lift-up effect arising in classical Navier Stokes case and the wave-shear coupled lift-up effects are suppressed, and there holds the following bound\[
\left\|\left(u,\theta\right)_{01}\right\|_{L^2}\lesssim C_{\widetilde{\gamma}}\left\|\left(u,\theta\right)_{01}(0)\right\|_{L^2}
\]where $C_{\widetilde{\gamma}}=\sqrt{\frac{\widetilde{\gamma}}{\widetilde{\gamma}-1}}+\sqrt{\frac{\widetilde{\gamma}-1}{\widetilde{\gamma}}}$. The dispersive estimate for $\widehat{u}_{2,01}$ and $\widehat{u}_{3,01}$ when $\sigma\neq1$ and $\sigma>0$ was proved in \cite{huang2026stabilitythreshold3dboussinesq}, we therefore omit the details here. We also note that if the initial PV vanishes, then the time-independent terms in \eqref{explicitsolution000} vanish, leaving $\widehat{u}_{1,01}$ and $\widehat{\theta}_{01}$ purely dispersive.

When $\widetilde{\gamma}\in(0,1)$, we have $\sigma<0$, and the dynamics are significantly different. Indeed, near the regime where $\eta\approx0$, one can verify that any $\widetilde{\gamma}\in(0,1)$ can lead to spectral instability such that\[
\left\|\left(u,\theta\right)_{01}\right\|_{L^2}\gtrsim e^{\Lambda t}\left\|\left(u,\theta\right)_{01}(0)\right\|_{L^2}
\]for some positive number $\Lambda>0$. 
\subsection*{4.3. Double zero modes} The dynamics of double zero modes are very simple. Due to the incompressibility and the assumption that $u_{2,00}\in L^2_y$, we have $u_{2,00}\equiv0$. The linearized system reduce to\[\left\{
\begin{aligned}
&\partial_tu_{1,00}=0,\\
&\partial_t u_{3,00}=\sqrt{\alpha}\theta_{00},\\
    &\partial_t\theta_{00}=-\sqrt{\alpha} u_{3,00}.
\end{aligned}
\right.\]We see that\[
||(u_1,u_3,\theta)_{00}||_{L^2}\lesssim||(u_1,u_3,\theta)_{00}(0)||_{L^2}.
\]
\section{Failure of a naive PV–wave transform}\label{section5}
 The PV-Wave decomposition method can be a powerful tool for analyzing system with scaling-compatible symmetrizor. However, generalizing this method to non-symmetric systems or systems without a scaling-compatible symmetrizer may not be straightforward. Let us consider the following example: setting $\widetilde{\gamma}=0$ in the Boussinesq system \eqref{2342342}, the $(U_3,\Theta,W_3)$ system satisfies
\begin{flalign*}
\left\{
\begin{aligned}
&\partial_t  \widehat{U}_3=\frac{f'l^2}{pf}\widehat{U}_3+\sqrt{\alpha}\frac{f}{p}\widehat{\Theta}-\frac{2ik^2l}{pf}\widehat{W}_3,\\
&\partial_t\widehat{W}_3=-il\widehat{U}_3,\\
    &\partial_t\widehat{\Theta}=-\sqrt{\alpha}\widehat{U}_3.
\end{aligned}
\right.
\end{flalign*}This system does not have scaling-compatible symmetrizer. Since the correct scaling is $\widehat{W}_3\sim 1$ and $\frac{2ik^2l}{pf}\sim\frac{1}{t^4}$ (see \cite{yao2026stabilitystratifiedcouette}), the term $-\frac{2ik^2l}{pf}\widehat{W}_3$ can only be treated as an error term. If we follow the steps introduced in the previous sections, we can identify the wave frequency $\omega:=\sqrt{\alpha}\sqrt{\frac{f}{p}}$. Besides, we apply the symmetrization for $\Theta$ and $U_3$:\[
\widehat{\Theta}=A,\qquad \widehat{U}_3=\sqrt{\frac{f}{p}}B.
\]We can then identify the PV ($K$) and Wave ($Y$) unknowns:\[
K=\widehat{W}_3-\frac{il}{\sqrt{\alpha}}A,\qquad Y=J(t)\left(-\frac{il}{\sqrt{\alpha}}\widehat{W}_3+A\right).
\]To symmetrize the $(B,Y)$ system, we need to pick $
J=\frac{\alpha}{\alpha+l^2}.$ Consequently, the $(B,Y,K)$ system is formulated as\begin{flalign*}
\left\{
\begin{aligned}
&\partial_t K=0,\\
&\partial_tY=-\sqrt{\alpha}\sqrt{\frac{f}{p}}B,\\
    &\partial_tB=\frac{f'l^2}{2fp}B+\sqrt{\alpha}\sqrt{\frac{f}{p}}Y+\frac{il\alpha}{\alpha+l^2}\sqrt{\frac{f}{p}} K+\frac{2k^2l^2}{\sqrt{\alpha p}f^\frac{3}{2}}Y-\frac{2ik^2l}{\sqrt{p}f^\frac{3}{2}}\frac{\alpha}{\alpha+l^2}K.
\end{aligned}
\right.
\end{flalign*}Due to the presence of the term $\frac{il\alpha}{\alpha+l^2}\sqrt{\frac{f}{p}}K$, whose coefficient is non-time-integrable. Hence this transformation does not seem to provide any simplification. A method to address this system is to reformulate the system $(u_3,\theta)$ directly using the linearized PV $K$. Both the linear and nonlinear stability can be proved, see \cite{yao2026stabilitystratifiedcouette} for details.

\appendix
\section{Appendix}
\subsection*{Magnetorotational instability}\label{MRI}

Consider the system \eqref{RBM} linearized around the following steady state\[
v^*=(y,0,0),\qquad B^*=(0,0,\beta),\qquad \varrho^*=\alpha z.
\]We write the perturbations $v=v^*+u$, $B=B^*+b$ and $\varrho=\varrho^*+\sqrt{\alpha}\theta$ and consider the $x$-averaged part: \begin{flalign}
\left\{
\begin{aligned}
& \partial_t u+\begin{pmatrix}
    (1-\widetilde{\gamma})u_2\\\widetilde{\gamma}u_1\\0
\end{pmatrix}+\bm{\nabla}_{y,z}\Delta^{-1}_{y,z}\left(\sqrt{\alpha}\partial_z\theta-\widetilde{\gamma}\partial_y u_1\right)=\beta\partial_zb+\begin{pmatrix}
    0\\0\\\sqrt{\alpha}\theta
\end{pmatrix},\\
& \partial_t\theta+\sqrt{\alpha}u_3=0,\\
& \partial_t b=\beta\partial_zu+\begin{pmatrix}
    b_2\\0\\0
\end{pmatrix},\\
& \partial_y u_2+\partial_zu_3=0,\qquad \partial_y b_2+\partial_z b_3=0.
\end{aligned}
\right.
\label{xLRSM}
\end{flalign}We consider the perturbations defined on $(y,z)\in\mathbb{R}\times\mathbb{T}$ and further restrict the solutions without $z$-average. Taking the Fourier transform, we can derive an ODE for $u_2$:\[
\partial^4_t\widehat{u}_2+\left(2\beta^2l^2+\frac{\alpha\eta^2+\widetilde{\gamma}(\widetilde{\gamma}-1)l^2}{\eta^2+l^2}\right)\partial^2_t\widehat{u}_2+\beta^2l^2\left(\beta^2l^2+\frac{\alpha\eta^2-\widetilde{\gamma}l^2}{\eta^2+l^2}\right)\widehat{u}_2=0.
\]Let $\widehat{u}_2\sim e^{\mu t}$, $\mu$ solves the following equation\begin{equation}
\mu^4+\left(2\beta^2l^2+\frac{\alpha\eta^2+\widetilde{\gamma}(\widetilde{\gamma}-1)l^2}{\eta^2+l^2}\right)\mu^2+\beta^2l^2\left(\beta^2l^2+\frac{\alpha\eta^2-\widetilde{\gamma}l^2}{\eta^2+l^2}\right)=0.\label{eigen}
\end{equation}Now we determine when the equation has a root with positive real part. Suppose first that some roots are imaginary of the form $\mu^2=re^{i\theta}$ for some $r\geq0$, then we have $\mu=\sqrt{r}e^{i\frac{\theta}{2}+ik\pi}$ for $k=0,1$. Thus, the only possibility for stability is that $\mu^2<0$, i.e. the roots are negative (or non-positive to include the neutral stability). Consider the quadratic equation (set $\lambda=\mu^2$)\[
\lambda^2+a\lambda+b=0,
\]the equation has two distinct negative roots if and only if\[
a>0,\qquad b>0,\qquad a^2-4b>0.
\]It can be verified that $a^2-4b>0$ always holds true. For $b>0$, by writing $\eta=l\tan\theta$, we deduce that\[
b>0,\qquad \forall\theta\in[0,2\pi]\iff \widetilde{
\gamma}<\beta^2.
\]Similarly, for $a>0$, this would require (when $1-8\beta^2>0$)\[
-\widetilde{\gamma}(\widetilde{\gamma}-1)<2\beta^2\iff \left|\widetilde{\gamma}-\frac{1}{2}\right|>\frac{\sqrt{1-8\beta^2}}{2}.
\]When $1-8\beta^2<0$, $a>0$. Since $\frac{1-\sqrt{1-8\beta^2}}{2}>\beta^2$, we see that if $b>0$, then $a>0$. Therefore, we conclude that the solution is only stable if $\widetilde{\gamma}<\beta^2$.

\section*{Statements and Declarations}

\noindent\textbf{Competing interests.}
The author affirms that he has no conflict of interest regarding this work.\vspace{5pt}

\noindent\textbf{Data and code availability.}
The Python code used to generate the illustrative numerical figures is available from the corresponding author upon reasonable request. No separate datasets were generated.
\printbibliography
\end{document}